\documentclass[letterpaper]{article}

\usepackage{amssymb}
\usepackage{amsfonts}
\usepackage{amsmath}
\usepackage{amscd}
\usepackage{amsmath}
\usepackage{amsxtra}
\usepackage{amsthm}

\theoremstyle{plain}
\newtheorem{theorem}{Theorem}[section]

\newtheorem{corollary}[theorem]{Corollary}

\newtheorem{lemma}[theorem]{Lemma}

\newtheorem{proposition}[theorem]{Proposition}

\theoremstyle{definition}
\newtheorem{definition}[theorem]{Definition}
\newtheorem{example}[theorem]{Example}

\newtheorem{remark}[theorem]{Remark}

\input{tcilatex}
\numberwithin{equation}{section}

\begin{document}

\title{Wavelets on Fractals\thanks{%
Work supported by the National Science Foundation.\newline
} \thanks{%
AMS Subject Classification: [2000] 41A15, 42A16, 42A65, 42C40, 43A65, 46L60,
47D25, 46L45. \newline
\qquad Keywords: Hausdorff measure, Cantor sets, iterated function systems
(IFS), fractal, wavelets, Hilbert space, unitary operators, orthonormal
basis (ONB), spectrum, transfer operator, cascade approximation, scaling,
translation.\newline
}}
\author{Dorin E. Dutkay\thanks{%
ddutkay@math.uiowa.edu\newline
} and Palle E. T. Jorgensen\thanks{%
jorgen@math.uiowa.edu\newline
} \\
The University of Iowa} \maketitle \begin{center}{\it Dedicated to
the memory of Gert Kj\ae rgaard Pedersen}\end{center}
\begin{abstract}
We show that there are Hilbert spaces constructed from the Hausdorff
measures $\mathcal{H}^{s}$ on the real line $\mathbb{R}$ with $0<s<1$ which
admit multiresolution wavelets. For the case of the middle-third Cantor set $%
\mathbf{C}\subset \lbrack 0,1]$, the Hilbert space is a separable
subspace of $L^{2}(\mathbb{R},(dx)^{s})$ where $s=\log _{3}(2)$.
While we develop the general theory of multi-resolutions in
fractal Hilbert spaces, the emphasis
is on the case of scale $3$ which covers the traditional Cantor set $\mathbf{%
C}$. Introducing
\begin{equation*}
\psi _{1}(x)=\sqrt{2}\chi _{\mathbf{C}}(3x-1)
\end{equation*}%
and
\begin{equation*}
\psi _{2}(x)=\chi _{\mathbf{C}}(3x)-\chi _{\mathbf{C}}(3x-2)
\end{equation*}%
we first describe the subspace in $L^{2}(\mathbb{R},(dx)^{s})$ which has the
following family as an orthonormal basis (ONB):%
\begin{equation*}
\psi _{i,j,k}(x)=2^{\frac{j}{2}}\psi _{i}(3^{j}x-k)\text{,}
\end{equation*}%
where $i=1,2,j$,$~k\in \mathbb{Z}$.
\par
Since the affine iteration systems of Cantor type arise from a
certain algorithm in $\mathbb{R}^d$ which leaves gaps at each
step, our wavelet bases are in a sense gap-filling constructions.
\end{abstract}

\tableofcontents

\section{\label{Intro}Introduction}
The paper has three interrelated themes: (1) construction of
wavelet bases in separable Hilbert spaces built on affine fractals
and Hausdorff measure; (2) approximation of the corresponding
wavelet scaling functions, using the cascading approximation
algorithm; and (3) an associated spectral theoretic analysis of a
transfer operator, often called the Ruelle operator.

There are surprises when our results are compared to what is known
for the traditional multiresolution approach for
$L^{2}(\mathbb{R}^{d})$, and even when compared to known results
for special classes of affine fractals.

Some comments on (1) - (3): Due to earlier work by Jorgensen,
Pedersen, \cite{JoPe98} and Strichartz et al \cite{Strich00}, it
is known that a subclass of the affine fractals admits Fourier
duality. Affine fractals arise from the specification of an
expansive matrix, and a finite set of translations. The fractal
$X$ itself then arises from this data and an iteration 'in the
small' of the corresponding affine maps. Let $L=L(X)$ be the
associated iteration 'in the large'. We say that $(X,L)$ is a
Fourier duality, if an orthonormal basis on $X$ may be built from
the frequencies in $L$. While it is known that, if $X$ is the
middle third Cantor set, then there is no $L$ which makes a
duality pair, we show that nonetheless, every affine fractal
admits an orthonormal wavelet basis. In our discussion of
wavelets, we start with the middle third Cantor set; and we then
pass on to the general affine fractals.

As for the approximation issues in (2), we know that for $L^{2}(\mathbb{R}%
^{d})$, there is a rich family of wavelet filters which yield
cascade approximation. This family of filters is much more
restricted for the fractals: Our results for the affine fractals
even offer a certain dichotomy (Theorem \ref{XformTheo1}): If the
cascades do not converge in the Hilbert space, then the terms in
the cascading approximation sequence are typically orthogonal, and
thus very far from being convergent.

Our analysis of (1) - (2) is based on spectral theory of the
associated transfer operator, and we show in the second half of
the paper (starting with Section \ref{ZakTransform}) how this
spectral theory differs in the three cases, the standard
$L^{2}(\mathbb{R}^{d})$- wavelets, and the special duality
fractals versus the general class of affine fractals.

Our proofs depend on ideas from geometric measure theory, and from
earlier papers on harmonic analysis of affine fractals. While some
of this material is in the literature, it isn't available
precisely in the form we need it here. In any case, it is
difficult for readers to locate without first having a brief
overview. So we include a minimum amount of facts from the
literature for the benefit of readers. We hope thereby to bridge
the diverse fields, fractals, Hilbert space, wavelets,
approximation, and harmonic analysis.

We develop the theory of multiresolutions in the context of Hausdorff
measure of fractional dimension between 0 and 1. While our fractal wavelet
theory has points of similarity that it shares with the standard case of
Lebesgue measure on the line, there are also sharp contrasts. These are
stated in our main result, a dichotomy theorem. The first section is the
case of the middle-third Cantor set. This is followed by a review of the
essentials on Hausdorff measure. The remaining sections of the paper cover
multiresolutions in the general context of affine iterated function systems.

It is well known that the Hilbert spaces $L^{2}(\mathbb{R})$ has a rich
family of orthonormal bases of the following form:%
\begin{equation*}
\psi _{j,k}(x)=2^{j/2}\psi (2^{j}x-k),\qquad j,k\in \mathbb{Z},
\end{equation*}%
where $\psi $ is a single function $\in L^{2}(\mathbb{R})$, with
\begin{equation*}
\left\Vert \psi \right\Vert _{2}=\left( \int_{\mathbb{R}}\left\vert \psi
(x)\right\vert ^{2}dx\right) ^{1/2}=1\text{,}
\end{equation*}%
and the integration refers to the usual Lebesgue measure on $\mathbb{R}$.
Take for example%
\begin{equation}
\psi (x)=\chi _{I}(2x)-\chi _{I}(2x-1)  \label{IntroEq1}
\end{equation}%
where $I=[0,1]$ is the unit interval.

Clearly $I$ satisfies%
\begin{equation*}
2I=I\cup (I+1)\text{.}
\end{equation*}%
The Cantor subset $\mathbf{C}\subset I$ satisfies%
\begin{equation}
3\mathbf{C}=\mathbf{C}\cup (\mathbf{C}+2)  \label{IntroEq2}
\end{equation}%
and its indicator function $\varphi _{\mathbf{C}}:=\chi _{\mathbf{C}}$
satisfies%
\begin{equation}
\varphi _{\mathbf{C}}(\frac{x}{3})=\varphi _{\mathbf{C}}(x)+\varphi _{%
\mathbf{C}}(x-2)\text{.}  \label{IntroEq3}
\end{equation}

Since both constructions, the first one for the Lebesgue measure, and the
second one for the Hausdorff version $(dx)^{s}$, arise from scaling and
subdivision, it seems reasonable to expect multiresolution wavelets also in
Hilbert spaces constructed on the scaled Hausdorff measures $\mathcal{H}^{s}$
which are basic for the kind of iterated function systems which give Cantor
constructions built on scaling and translations by lattices. We show this to
be the case, but there are still striking differences between the two
settings, and we spell out some of them after first developing the theory in
the case of the middle-third Cantor construction.

While there are other wavelet approaches to fractals in the
literature, for example \cite{Jon98}, \cite{Jon04}, and
\cite{Strich97}, there is in fact no overlap with this work, since
the previous papers deal with wavelets on the fractal itself,
while the present paper deals with wavelets on an enlarged fractal
(actually a fractal measure), allowing a structure closer to a
standard multiresolution analysis (MRA).

The practical applications are to fractals arising in physics and
in symbolic dynamical systems from theoretical computer science,
see e.g., \cite{Berry} \cite{Sadun}, \cite{Takens}. There is
already a considerable body of work on harmonic analysis on
fractals, see for example \cite{Strich03}, \cite{HuStrich}, \cite{GibRajStr}%
, \cite{Strich00}, and \cite{JoPe98}. Much of it is based on subdivision
techniques, and algorithms which use cascade constructions, but so far we
have not seen direct wavelet algorithms and wavelet analysis for fractals.

In section 2, we recall some facts about Hausdorff measure $%
\mathcal{H}^{s}$, Hausdorff dimension, and Hausdorff distance. They will be
needed in the Hilbert space we build on $\mathcal{H}^{s}$. It is a natural
separable subspace of the full $\mathcal{H}^{s}$-Hilbert space, and it is
built up from the algebra of $\mathbb{Z}$-translations (additive), and $N$%
-adic scaling (multiplicative), where $N$ is fixed. We then turn to the
cascade approximation for the scaling function $\varphi $ defined by the
usual $1/N$ subdivision$.$ We prove a theorem for the case $0<s<1$ which
stands in sharp contrast to the traditional and more familiar case $s=1$ of
Daubechies et. al.; i.e., the case of the Hilbert space $L^{2}(\mathbb{R)}$
based on Lebesgue measure $dx$ on $\mathbb{R}$: The scaling equation is then
\begin{equation}
\varphi (x)=\sqrt{N}\dsum\limits_{k\in \mathbb{Z}}a_{k}\varphi (Nx-k)
\label{IntroEq3a}
\end{equation}%
with the masking coefficients $a_{k}$ satisfying the usual two axioms
\begin{equation}
\dsum\limits_{k\in \mathbb{Z}}a_{k}=\sqrt{N}\text{, and }\dsum\limits_{k\in
\mathbb{Z}}\bar{a}_{k}a_{k+N\ell }=\delta _{\ell ,0}\text{, }\ell \in
\mathbb{Z}\text{.}  \label{IntroEq3b}
\end{equation}%
Motivated by the expression on the right hand side in (\ref{IntroEq3a}), we
define the wavelet subdivision operator $M$ by
\begin{equation}
(Mf)(x):=\sqrt{N}\dsum\limits_{k\in \mathbb{Z}}a_{k}f(Nx-k)\text{,\qquad }%
f\in L^{2}(\mathbb{R)}\text{;}  \label{IntroEq3c}
\end{equation}%
and note that its properties depend on the specifications in (\ref{IntroEq3b}%
).

Simple conditions are known for when the limit
\begin{equation}
\lim_{n\rightarrow \infty }M^{n}\chi _{I}=\varphi  \label{IntroEq3d}
\end{equation}%
exists in $L^{2}(\mathbb{R)}$, see \cite{Dau92}, chapter 5. Then $\varphi $ (when it exists) solves (\ref%
{IntroEq3a}), and there is an easy formula for building functions
$\psi _{1},\ldots ,\psi _{N-1}$ in $L^{2}(\mathbb{R)}$ from
$\varphi $ such that
\begin{equation}
\left\{ N^{\frac{k}{2}}\psi _{i}(N^{k}x-\ell )\mid 1\leq i<N,k,\ell \in
\mathbb{Z}\right\}  \label{IntroEq3e}
\end{equation}%
is an orthonormal basis (ONB) in $L^{2}(\mathbb{R)}$. If $N=2$, a formula
for $\psi $ is
\begin{equation}
\psi (x)=\sqrt{2}\dsum\limits_{k\in
\mathbb{Z}}(-1)^{k}\bar{a}_{1-k}\varphi (2x-k)\text{.}
\label{IntroEq3f}
\end{equation}%
In general when $N\geq 2$, the functions $\psi _{1},\ldots ,\psi
_{N-1}$ may result from the solution to a simple matrix completion
problem; see \cite{PR03}, \cite{BrJo01a} and \cite{BrJo02} for details. In the case of Hausdorff measure $%
\mathcal{H}^{s}$ ($0<s<1$, $s$ depending on the scaling number
$N$), the analogous matrix completion is still fairly simple. A
main question (non-trivial) is now that of solving the analogue of
(\ref{IntroEq3a}), but in the $\mathcal{H}^{s}$-Hilbert space. The
biggest differences concern the changes in (\ref{IntroEq3b}) and
(\ref{IntroEq3d}) when $0<s<1$. It turns out in the fractal cases
that there are then many fewer admissible solutions than those
suggested by analogy with (\ref{IntroEq3b}). We summarize the
situation in Sections 4--6, where our main result takes the form
of a dichotomy theorem; the solutions to the
$\mathcal{H}^{s}$-convergence question are isolated within a
larger family of masking coefficients analogous to
(\ref{IntroEq3b}). There is further a new set of orthogonality
conditions entering the analysis when $0<s<1$, which are not
present in the more familiar case of $s=1$.

We interpret the wavelet filters as functions $m_{0}$ on the torus $\mathbb{T%
}$. If the scaling number $N$ is given, following \cite{Jorgen01},
we introduce the wavelet-transfer operator
\begin{equation}
\left( R_{m_{0}}f\right) (z):=\frac{1}{N}\dsum\limits_{w^{N}=z}\left\vert
m_{0}(w)\right\vert ^{2}f(w)\text{,\qquad for }f\in C(\mathbb{T)}\text{, and
}z\in \mathbb{T}\text{.}  \label{IntroEq3g}
\end{equation}

Our dichotomy for wavelets will be explained in terms of the
spectral properties of $R_{m_{0}}$, also called the Ruelle
operator. For simplicity, we introduce the normalization
$R_{m_{0}}(\hat{1})=\hat{1}$, where $\hat 1$ denotes the constant
function $1$ on $\mathbb{T}$. A probability measure $\nu $ on
$\mathbb{T}$ is said to be \textit{invariant} if $\nu
R_{m_{0}}=\nu $. Equivalently,
\begin{equation*}
\dint\limits_{\mathbb{T}}R_{m_{0}}(f)d\nu =\dint\limits_{\mathbb{T}}fd\nu
\text{,\qquad for all }f\in C\left( \mathbb{T}\right) \text{,}
\end{equation*}%
or
\begin{equation}
\nu (R_{m_{0}}(f))=\nu \left( f\right) \text{.}  \label{IntroEq3h}
\end{equation}%
We introduce a notion of $\left( m_{0},N\right) $-cycles for (\ref{IntroEq3g}%
) which explains the solutions $\nu \in M_{1}\left( \mathbb{T}\right) $ to (%
\ref{IntroEq3h}). In our setting, the dichotomy boils down to two cases for $%
\nu $:\smallskip

\noindent (i) $\nu =\delta _{1}$ (the Dirac mass at $z=1$), or

\noindent (ii) some $\nu $ is a singular measure on $\mathbb{T}$ with full
support.\medskip

\noindent In the first case (i), the Hilbert space is
$L^{2}(\mathbb{R)}$; i.e., that of the standard wavelets; and in
the second case (ii), the Hilbert space is built from Hausdorff
measure $\mathcal{H}^{s}$, $0<s<1$.

We begin the discussion with $N=3$ and $s=\log _{3}(2)$.

\begin{definition}
\label{IntroDef1}We define $\mathcal{R}$ to be the set of all real
numbers that have a base $3$ expansion containing only finitely
many ones. It is an inflated version of $\mathbf{C}$:
\end{definition}

$\mathcal{R}$:$=\{\sum_{k=-m}^{\infty }a_{k}3^{-k}\mid m\in \mathbb{Z}%
,a_{k}\in \{0,1,2\}$ for all $k\in \mathbb{Z},a_{k}\neq 1$ for all but
finitely many indices $k\}$.

For the fractal cases, $0<s<1$, the factor $\sqrt{N}$ in equations (\ref%
{IntroEq3a})--(\ref{IntroEq3c}) will be different, see details below.
Similarly the factor $N^{\frac{k}{2}}$ in (\ref{IntroEq3e}) changes: With
scaling number $N$ and with $p$ subdivisions, the ONB corresponding to (\ref%
{IntroEq3e}) in $L^{2}(\mathcal{R},\mathcal{H}^{s})$, $s=\log _{N}(p)$ is $%
\{p^{\frac{k}{2}}\psi _{i}(N^{k}x-\ell )\}$. However, the geometric
properties of the cascade approximation to the scaling function change
completely as the Hausdorff dimension moves from $s=1$ to the open interval $%
0<s<1$. This will be spelled out in the last four sections of the paper.

\section{\label{HausMeasure}The Hausdorff Measure}
Returning to the middle-third Cantor set $\mathbf{C}=\mathbf{C}_3$; i.e., $N=3$ and $p=2$%
, here are some elementary properties of $\mathcal{R}$:

\begin{proposition}
\label{IntroProp1}\textup{(}The middle-third Cantor set.\textup{)} The set $%
\mathcal{R}$ has the following properties:

\noindent \textup{(}i\textup{)} Invariance under triadic translation:%
\begin{equation*}
\mathcal{R}+\frac{k}{3^{n}}=\mathcal{R},\qquad (k,n\in \mathbb{Z})\text{.}
\end{equation*}

\noindent \textup{(}ii\textup{)} Invariance under dilation by $3$:%
\begin{equation*}
3^{n}\mathcal{R}=\mathcal{R},\qquad (n\in \mathbb{Z}).
\end{equation*}

\noindent \textup{(}iii\textup{)} The middle-third Cantor set $\mathbf{C}$
is contained in $\mathcal{R}$ and moreover it covers $\mathcal{R}$ by
translations and dilations:%
\begin{equation}
\mathcal{R}=\dbigcup_{n\in\mathbb{Z}}\dbigcup\limits_{k\in \mathbb{Z}}3^{-n}(%
\mathbf{C}+k). \label{IntroEq4}
\end{equation}
\end{proposition}

\begin{proof}
\noindent (i) Any triadic number $t=\frac{k_{0}}{3^{n_{0}}}$ with $%
k_{0},n_{0}\in \mathbb{Z}$, $k_{0}\geq 0$, has a finite expansion in base $3$%
:%
\begin{equation*}
t=\sum_{k=-m}^{m}t_{k}3^{-k}\text{,}\qquad t_{k}\in \left\{ 0,1,2\right\}
\text{.}
\end{equation*}

Take $x\in \mathcal{R}$. Then $x$ has a finite number of ones in its
expansion so the same affirmation will be true for $x+t$.

\noindent (ii) is clear: multiplication by $3$ means a shift in the base $3$
expansion.

\noindent (iii) Since%
\begin{equation}
\mathbf{C}=\left\{ \sum_{k=1}^{\infty }a_{k}3^{-k}\mid a_{k}\in
\{0,2\}\right\} \text{,}  \label{IntroEq5}
\end{equation}%
it is obvious that $\mathbf{C}\subset \mathcal{R}$.

The inclusion \textquotedblleft $\supset $\textquotedblright\ follows from
(i) and (ii). Now take
\begin{equation*}
x\in \mathcal{R}\text{,\quad }x=\sum_{k=-m}^{\infty }a_{k}3^{-k}\text{,}
\end{equation*}%
only finitely many $a_{k}$ are equal to $1$.

Let $\ell _{0}$ be the last index for which $a_{l_{0}}=1$, and take $k_{0}$:~%
$=\max\{n_{0},\ell _{0}\}$.

Then
\begin{equation*}
x\in 3^{-k_{0}}\left( \mathbf{C}+\sum_{k=-m}^{k_{0}}a_{k}3^{-k+k_{0}}\right)
\end{equation*}%
and this shows that the other inclusion is also true.
\end{proof}

\begin{remark}
$\mathcal{R}$ has Lebesgue measure $0.$ Indeed, this follows from
proposition 2 \textup{(}iii\textup{)}, because $\mathbf{C}$ has Lebesgue
measure $0$, and so do all the sets $3^{-n}(\mathbf{C}+k)$ with $n,k\in
\mathbb{Z}$.
\end{remark}

Even though some of the properties that we need for the Hausdorff
measure, for fractals, and for iterated function systems (IFS) are
known, we found that the material is wildly scattered throughout
the literature; and to increase readability we have included some
highpoints from these areas. This should also help bring out the
contrast between the traditional MRA-analysis, and the present
more stochastic approach.

Next we define a measure on $\mathcal{R}$. It is the restriction of the
Hausdorff measure $\mathcal{H}^{s}$ with $s=\log _{3}(2)$ to $\mathcal{R}$.

We recall some background on the Hausdorff measures from \cite{Fal}:

For a subset $E$ of $\mathbb{R}$, $s>0$, and $\delta >0$, define%
\begin{equation*}
\mathcal{H}_{\delta }^{s}(E):=\inf \left\{ \sum_{i=1}^{\infty }\left\vert
U_{i}\right\vert ^{s}\mid \dbigcup\limits_{i=1}^{\infty }U_{i}\supset
E,\left\vert U_{i}\right\vert <\delta \right\}
\end{equation*}%
where $\left\vert U\right\vert =\sup \left\{ \left\vert x-y\right\vert \mid
x,y\in U\right\} $ (the diameter of $U$).

It is known that $\mathcal{H}_{\delta }^{s}$ is an outer measure on $%
\mathbb{R}
$.

Define%
\begin{equation}
\mathcal{H}^{s}(E)=\lim_{\delta \rightarrow 0}\mathcal{H}_{\delta
}^{s}(E)=\sup_{\delta >0}\mathcal{H}_{\delta }^{s}(E)\text{.}
\label{HausEq1}
\end{equation}%
Then verify that $\mathcal{H}^{s}$ is also an outer measure. By the
Caratheodory construction \cite{Fal}, if we restrict $\mathcal{H}^{s}$ to
the $\sigma $-field of $\mathcal{H}^{s}$-measurable sets, we get a measure
called the Hausdorff measure.

\begin{proposition}
\label{HausMeasProp1}\emph{:}

\noindent \textup{(}i\textup{)} All Borel sets are measurable.\emph{\ }

\noindent \textup{(}ii\textup{)} \textup{[}Inner regularity\textup{]} Any $%
\mathcal{H}^{s}$-measurable set of finite $\mathcal{H}^{s}$-measure contains
a $F_{\sigma }$-set of equal $\mathcal{H}^{s}$-measure.

\noindent \textup{(}iii\textup{)} If $E\subset \mathcal{R}$\ then there is a
$G_{\delta }$-set containing $E$\ and of the same $\mathcal{H}^{s}$-measure.

\noindent \textup{(}iv\textup{)} For $s<1$\ and $G$\ open, $\mathcal{H}%
^{s}(G)=\infty $\ \textup{(}the measure is not regular from above\textup{)}.

\noindent \textup{(}v\textup{)} Translation invariance: For any $\mathcal{H}%
^{s}$-measurable set $E$, and any $t\in \mathcal{R}$, $E+t$\ is $\mathcal{H}%
^{s}$-measurable, and\textit{\ }%
\begin{equation}
\mathcal{H}^{s}(E)=\mathcal{H}^{s}(E+t).  \label{HausEq1a}
\end{equation}

\noindent \textup{(}vi\textup{)} For any $\mathcal{H}^{s}$-measurable set $E$%
\ and any $c>0$, $cE$\ is $\mathcal{H}^{s}$-measurable and\textit{\ }%
\begin{equation}
\mathcal{H}^{s}(cE)=c^{s}\mathcal{H}^{s}(E)\mathit{.}  \label{HausEq1b}
\end{equation}
\end{proposition}

Consider now $\mathcal{H}^{s}$ with $s=\log _{3}2$ restricted to the $%
\mathcal{H}^{s}$-measurable subsets of $\mathcal{R}$. We will keep the
notation $\mathcal{H}^{s}$ for the restriction.

\begin{proposition}
\label{HausMeasProp2}\emph{:}

\noindent \textup{(}\textit{i}\textup{)} \textit{If }$E\subset \mathcal{R}$%
\textit{\ is an }$\mathcal{H}^{s}$\textit{-measurable set and }$t=\frac{l_{0}%
}{3^{p_{0}}}$\textit{\ is a triadic number, then }$E+t\subset R$\textit{\ is
}$\mathcal{H}^{s}$\textit{-measurable and}
\begin{equation*}
\mathcal{H}^{s}(E)=\mathcal{H}^{s}(E+t)\text{.}
\end{equation*}

\noindent \textup{(}\textit{ii}\textup{)} \textit{If }$E\subset \mathcal{R}$%
\textit{\ is an }$\mathcal{H}^{s}$\textit{-measurable set then }$3E\subset
\mathcal{R}$\textit{\ is }$\mathcal{H}^{s}$\textit{-measurable, and }%
\begin{equation}
\mathcal{H}^{s}(3E)=2\mathcal{H}^{s}(E)\mathit{.}  \label{HausEq1c}
\end{equation}

\noindent \textup{(}\textit{iii}\textup{)} \textit{If }$f\in L^{1}(\mathcal{R%
},\mathcal{H}^{s}$\textit{) then the function on }$\mathcal{R}$\textit{, }$%
x\rightarrow f(\frac{x}{3})$\textit{\ is also in }$L^{1}(\mathcal{R},%
\mathcal{H}^{s}$\textit{) and}%
\begin{equation*}
\dint\limits_{\mathcal{R}}f(x)d\mathcal{H}^{s}(x)=\frac{1}{2}\dint\limits_{%
\mathcal{R}}f(\frac{x}{3})d\mathcal{H}^{s}(x)\text{.}
\end{equation*}

\noindent \textup{(}\textit{iv}\textup{)} $\mathcal{H}^{s}(\mathbf{C})=1$,
where $\mathbf{C}$ is the middle-third Cantor set.
\end{proposition}

\begin{proof}
\noindent (i) and (ii) are direct consequences of propositions \ref%
{HausMeasProp1} and \ref{HausMeasProp2}.

\noindent (iii) follows from (ii) when $f$ is a characteristic function of
an $\mathcal{H}^{s}$-measurable set. Then, for arbitrary $f$, the formula
can be obtained by approximations by simple functions.

\noindent (iv) See \cite{Fal} theorem 1.14.
\end{proof}

\begin{remark}
The measure $\mathcal{H}^{s}$ on $\mathcal{R}$ is still non-regular from
above. All open sets in $\mathcal{R}$ still have infinite measure.
\end{remark}

To see this, we show that $\mathcal{H}^{s}(I)=\infty $, where $I=(0,1)\cap
\mathcal{R}$.

Indeed $I\supset \mathbf{C}$, so $\mathcal{H}^{s}(I)\geq 1$.

Also, observe that $3I=I\cup (I+1)\cup (I+2)$ disjoint union (we neglect
some points that have $\mathcal{H}^{s}$-measure $0$.).

Therefore, with propositions \ref{HausMeasProp2} (i) and (ii), we obtain%
\begin{equation}
2\mathcal{H}^{s}(I)=3\mathcal{H}^{s}(I)  \label{HausEq2}
\end{equation}%
so $\mathcal{H}^{s}(I)$ is either $0$ or $\infty $. $0$ cannot be from the
previous argument, hence it must be $\infty $.

By scalings and translations, it can be proved that $\mathcal{H}%
^{s}((a,b)\cap \mathcal{R})=\infty $ for any interval $(a,b)$.

Since all open subsets of $\mathcal{R}$ have measure $\infty $ it follows
that no non-zero continuous function on $\mathbb{R}$ is integrable! (Just
take $f^{-1}((a,b))$ for some interval that doesn't contain $0$ and
intersects the range.)

\begin{definition}
\label{HausMeasDef1}We denote by $H$ the Hilbert space%
\begin{equation*}
H:=L^{2}(\mathcal{R},\mathcal{H}^{s})\text{.}
\end{equation*}%
The linear operator $T$ on $H$ defined by%
\begin{equation}
Tf(x)=f(x-1)\text{,\qquad }(f\in H,x\in \mathcal{R})\text{,}  \label{HausEq3}
\end{equation}%
is called the translation operator. The linear operator $U$ on $H$ defined
by
\begin{equation}
Uf(x)=\frac{1}{\sqrt{2}}f\left( \frac{x}{3}\right) \text{,\qquad }(f\in
H,x\in \mathcal{R})  \label{HausEq4}
\end{equation}%
is called the dilation operator.
\end{definition}

From Proposition \ref{HausMeasProp2}, and using some simple
computations, we obtain the following proposition:

\begin{proposition}
\label{HausMeasProp3}\emph{:}

\noindent \textup{(}i\textup{)} $T$\ and $U$\ are unitary operators.

\noindent \textup{(}ii\textup{)} $UTU^{-1}=T^{3}$.
\end{proposition}

Denote by $\varphi =\chi _{\mathbf{C}}$, the characteristic function of the
Cantor set $C$. We prove that $\varphi $\ satisfies all the properties of a
scaling vector.

\begin{proposition}
\label{HausMeasProp4}The following hold\emph{:}

\noindent \textup{(}i\textup{)} \textup{[}The scaling equation\textup{]}\ $%
U\varphi =\frac{1}{\sqrt{2}}(\varphi +T^{2}\varphi )$.

\noindent \textup{(}ii\textup{)} \textup{[}Orthogonality of the translates%
\textup{]}\ $\left\langle T^{k}\varphi \mid \varphi \right\rangle =\delta
_{k},(k\in \mathbb{Z})$.

\noindent \textup{(}iii\textup{)} \textup{[}Cyclicity\textup{]}\ $\overline{%
\func{span}}\left\{ U^{n}T^{k}\varphi \mid n\in \mathbb{Z},k\in \mathbb{Z}%
\right\} =H$.
\end{proposition}

\begin{proof}
\par
\noindent (i) $U\varphi =\frac{1}{\sqrt{2}}\chi _{3\mathbf{C}}$, $%
T^{2}\varphi =\chi _{\mathbf{C}+2}$, but $3\mathbf{C}=\mathbf{C}\dbigcup (%
\mathbf{C}+2)$, so (i) follows.\newline
\noindent (ii) $T^{k}\varphi =\chi _{\mathbf{C}+k}$ so $T^{k}\varphi $ and $%
\varphi $ are disjointly supported for $k\neq 0$. For $k=0$, $\left\langle
\varphi \mid \varphi \right\rangle =\int_{\mathcal{R}}\chi _{\mathbf{C}}d%
\mathcal{H}^{s}=1$, by proposition \ref{HausMeasProp4} (iv).\newline
\noindent (iii) First take $E\subset R$ measurable and with $\mathcal{H}%
^{s}(E)<\infty $. We want to approximate $\chi _{E}$ by linear combinations
of functions of the form $U^{n}T^{k}\varphi $.

Note also that \begin{equation}\label{eq2_7_1} U^{n}T^{k}\varphi
=\chi _{3^{n}(\mathbf{C}+k)},\quad (n,k\in \mathbb{Z}).
\end{equation}
With proposition \ref%
{HausMeasProp2} \textup{(}iii\textup{)}, and partitioning $E$\ if necessary,
we may assume that $E$\ is contained in a set of the form $3^{n^0}(%
\mathbf{C}+k_{0})$. Applying dilations and translations we may further
assume that $E\subset \mathbf{C}$
\end{proof}

Define%
\begin{equation}
\mathcal{V}\text{:}=\left\{ C_{n,a_{n},\ldots ,a_{1}}\text{:}=3^{-n}\mathbf{C%
}+\sum_{k=1}^{n}a_{k}3^{-k}\mid a_{k}\in \left\{ 0,2\right\} n\geq 1\right\}
\text{.}  \label{HausEq5}
\end{equation}%
This family $\mathcal{V}$ is a Vitali class for $E$; i.e., for each $x\in E$
and each $\delta >0$, there is a $U\in \mathcal{V}$ with $x\in U$, and $%
0<\left\vert U\right\vert \leq \delta $.

Indeed, we see that for all $n\geq 1$:%
\begin{equation*}
\dbigcup\limits_{a_1,...,a_n\in \{0,2\}}C_{n,a_{n},\ldots ,a_{1}}=%
\mathbf{C}\text{.}
\end{equation*}%
Also using proposition \ref{HausMeasProp4}

\noindent (\ref{HausEq5}) $\mathcal{H}^{s}($\textbf{$C$}$_{n,a_{n},\ldots
,a_{1}})=\mathcal{H}^{s}(3^{-n}\mathbf{C})=(3^{-n})^{s}\mathcal{H}^{s}(%
\mathbf{C})=2^{-n}=\left\vert C_{n,a_{n},\ldots ,a_{1}}\right\vert ^{s}$. We
conclude that $\mathcal{V}$ is indeed a Vitali class for any subset $E$ of $%
\mathbf{C}$.

Then, by Vitali's covering theorem (see \cite{Fal} theorem 1.10) for a fixed
$\varepsilon >0$, there exists a finite or countable disjoint sequence of
sets $\{U_{i}\}$ from $\mathcal{V}$ such that either $\sum \left\vert
U_{i}\right\vert ^{s}=\infty $, or $\mathcal{H}^{s}(E\backslash \dbigcup
U_{i})=0$, and also
\begin{equation*}
\mathcal{H}^{s}(E)<\sum_{i}\left\vert U_{i}\right\vert ^{s}+\varepsilon
\text{.}
\end{equation*}

Since the sets $U_{i}$ are mutually disjoint and contained in $\mathbf{C}$,
and using (i) in (\ref{HausEq5}), it follows that
\begin{equation}
\sum_{i}\left\vert U_{i}\right\vert ^{s}=\sum_{i}\mathcal{H}^{s}(U_{i})=%
\mathcal{H}^{s}(\dbigcup\limits_{i}U_{i})\leq \mathcal{H}^{s}(\mathbf{C})=1%
\text{.}  \label{HausEq6}
\end{equation}%
Therefore the other variant must be true:

with $U$:$=\dbigcup\limits_{i}U_{i}$%
\begin{equation*}
\mathcal{H}^{s}(E\backslash U)=0\text{.}
\end{equation*}

On the other hand
\begin{eqnarray*}
\mathcal{H}^{s}(U\backslash E) &=&\mathcal{H}^{s}(U)-\mathcal{H}^{s}(U\cap E)
\\
&=&\mathcal{H}^{s}(U)-(\mathcal{H}^{s}(E)-\mathcal{H}^{s}(E\backslash U)) \\
&=&\mathcal{H}^{s}(U)-\mathcal{H}^{s}(E)=\sum_{i}\mathcal{H}^{s}(U_{i})-%
\mathcal{H}^{s}(E) \\
&=&\sum_{i}\left\vert U_{i}\right\vert ^{s}-\mathcal{H}^{s}(E)\,<\varepsilon
.
\end{eqnarray*}

Also, observe that%
\begin{equation}
C_{n,a_{n},\ldots ,a_{1}}=3^{-n}\left( \mathbf{C}+\sum_{k=1}^{n}a_{k}3^{n-k}%
\right) \text{,}  \label{HausEq7}
\end{equation}%
so by (\ref{eq2_7_1}), $\chi _{C_{n,a_{n},\ldots
,a_{1}}}=U^{-n}T^{l}\varphi $ with $l=\sum_{k=1}^{n}a_{k}3^{n-k}$.

Therefore we see that all measurable sets $E\subset \mathcal{R}$ with $%
\mathcal{H}^{s}(E)<\infty $ are in the $\func{span}$ of $\{U^{n}T^{k}\varphi
\mid n,k\in \mathbb{Z}\}$. Since all integrable functions $f\in H$ can be
approximated by simple functions, it follows that
\begin{equation}
H=\overline{\func{span}}\left\{ U^{n}T^{k}\varphi \mid n,k\in \mathbb{Z}%
\right\} \text{.}  \label{HausEq8}
\end{equation}

\section{\label{IFS}Iterated Function Systems (IFS) and gap-filling wavelets}

The middle-third Cantor set $\mathbf{C}$ of Section 1 is a special
case of an Iterated Function System (IFS). It falls in the
subclass of the IFSs which are called \textit{affine.}
Specifically, let $d\in \mathbb{Z}_{+}$, and let $A$ be a $d\times
d$ matrix of $\mathbb{Z}$. Suppose that the eigenvalues $\lambda
_{i}$ of $A$ satisfy $\left\vert \lambda _{i}\right\vert >1$. Set
$N:=\left\vert \det A\right\vert .$ These matrices
are called \textit{expansive.} Then note that the quotient group $\mathbb{Z}%
^{d}\diagup A(\mathbb{Z}^{d})$ is of order $N$. A subset
$\mathcal{D}\subset \mathbb{Z}^{d}$ is said to represent the
$A$-residues if the natural quotient mapping
\begin{equation}
\gamma \text{:}~\mathbb{Z}^{d}\rightarrow \mathbb{Z}^{d}\diagup A(\mathbb{Z}%
^{d})  \label{IFSeq1}
\end{equation}%
restricts to a bijection $\gamma _{\mathcal{D}}$ of $\mathcal{D}$ onto $%
\mathbb{Z}^{d}\diagup A(\mathbb{Z}^{d})$. For example, if $d=1$,
and $A=3$,
then we may take either one of the two sets $\left\{ 0,1,2\right\} $ or $%
\left\{ 0,1,-1\right\} $ as $\mathcal{D}$. The IFSs which we shall
look at will be constructed from finite subsets
$\mathcal{S}\subset \mathbb{Z}^{d}$
which represent the $A$-residues for some given expansive matrix $A$. If $(A,%
\mathcal{S})$ is a pair with these properties, define the maps
\begin{equation}
\sigma _{s}(x):=A^{-1}(x+s)\text{, }s\in \mathcal{S}\text{, }x\in \mathbb{R}%
^{d}\text{.}  \label{IFSeq2}
\end{equation}%
Using a theorem of Hutchinson \cite{Hut81}, we conclude that there
is a
unique measure $\mu =\mu _{(A,\mathcal{S})}$ with compact support $\mathbf{%
C=C}_{(A,\mathcal{S})}$ on $\mathbb{R}^{d}$ such that
\begin{equation}
\mu =\frac{1}{\#(\mathcal{S})}\dsum\limits_{s\in \mathcal{S}}\mu
\circ \sigma _{s}^{-1}\text{,}  \label{IFSeq3}
\end{equation}%
or equivalently
\begin{equation}
\int f(x)d\mu (x)=\frac{1}{\#(\mathcal{S})}\dsum\limits_{s\in \mathcal{S}%
}\int f(\sigma _{s}(x))d\mu (x)\text{.}  \label{IFSeq4}
\end{equation}%
The quotient mapping
\begin{equation}
\gamma \text{:~}\mathbb{R}^{d}\rightarrow \mathbb{T}^{d}:=\mathbb{R}%
^{d}\diagup \mathbb{Z}^{d}  \label{IFSeq5}
\end{equation}%
restricts to map $\mathbf{C}$ bijectively onto a compact subset of $\mathbb{T%
}^{d}$. The Hausdorff dimension $h$ of $\mu $ and of the support
$\mathbf{C}$ is
\begin{equation*}
h=\frac{\log \#(\mathcal{S})}{\log N}\text{.}
\end{equation*}%
The system $(\mathbf{C},\mu)$ is called a Hutchinson pair, see
lemma \ref{IFSlem1}. \par If $d=1$, we will look at two examples:
(i) $(A,\mathcal{S})=(3,\{0,2\})$
which is the middle-third Cantor set $\mathbf{C}$ in Section 1, and (ii) $(A,%
\mathcal{S})=(4,\{0,2\})$ which is the corresponding construction,
but starting with a subdivision of the unit interval $I$ into 4
parts, and in each step of the iteration omitting the second and
the fourth quarter interval. As noted, then
\begin{equation}
h_{(i)}=\log _{3}(2)=\frac{\log 2}{\log 3}\text{,\qquad and }h_{(ii)}=\frac{1%
}{2}\text{;}  \label{IFSeq6}
\end{equation}%
for more details, see \cite{JoPe98}.

Since the arguments from proposition \ref{HausMeasProp2} and
\ref{HausMeasProp4} generalize, we will only sketch the general
statements of results for the affine IFSs, those based on pairs
$(A,\mathcal{S})$ in $\mathbb{R}^{d}$ where the matrix $A$ and the
subset $\mathcal{S}\subset \mathbb{Z}^{d}$ satisfy the
stated conditions. The number $h$ will be $h=\frac{\log (\#(\mathcal{S}))}{%
\log \left\vert \det A\right\vert }$; i.e., the Hausdorff
dimension of the measure $\mu ,$ and its support $\mathbf{C}$
which are determined from the given pair $(A,\mathcal{S})$. We
will then be working with the corresponding Hausdorff measure
$\mathcal{H}^{h}$, but now as a measure defined on subsets of
$\mathbb{R}^{d}$. The facts from Section 2 apply also to this more
general case in $\mathbb{R}^d$, for example, property (1.2) for
the middle-third Cantor set, now takes the following form
\begin{equation}
A\mathbf{C}=\dbigcup\limits_{s\in \mathcal{S}}(\mathbf{C}+s)
\label{IFSeq7}
\end{equation}%
where $A\mathbf{C}:=\{Ax\mid x\in \mathbf{C}\}$, and $\mathbf{C}%
+s:=\{x+s\mid x\in \mathbf{C}\}$, or equivalently
\begin{equation}
\mathbf{C}=\dbigcup\limits_{s\in \mathcal{S}}\sigma
_{s}(\mathbf{C}) \label{IFSeq8}
\end{equation}%
where $\sigma _{s}(\mathbf{C}):=\{\sigma _{s}(x)\mid x\in
\mathbf{C}\}$. The conditions on the pair $(A,\mathcal{S})$
guarantees that the sets in the union on the right-hand side in
(\ref{IFSeq7}) or in (\ref{IFSeq8}), are mutually non-overlapping.
This amounts to the so-called open-set-condition of Hutchinson
\cite{Hut81}. The set $\mathcal{R}$ which is defined in
Proposition \ref{IntroProp1} in the special case of the
middle-third Cantor set is now instead
\begin{equation}
\mathcal{R}=\dbigcup\limits_{n\geq 0}\dbigcup\limits_{k\in \mathbb{Z}%
^{d}}A^{-n}(\mathbf{C}+k)=\mathcal{R}=\dbigcup\limits_{n\in\mathbb{Z}}\dbigcup\limits_{k\in \mathbb{Z}%
^{d}}A^{-n}(\mathbf{C}+k)    \label{IFSeq9}
\end{equation}%
where $\mathbf{C}$ is the (unique) compact set determined by
(\ref{IFSeq8}),
of Hutchinson's theorem \cite{Hut81}. The properties of Proposition \ref%
{IntroProp1} carry over \textit{mutatis mutandis}, for example,
the argument from Section 2 shows that for every $k\in
\mathbb{Z}^{d}$ and every $n\in \mathbb{Z}$,
\begin{equation*}
\mathcal{R}+A^{-n}k=\mathcal{R},\quad \text{and }A^{n}\mathcal{R}=\mathcal{R}%
.
\end{equation*}%
The Hilbert space $H$ from Definition \ref{HausMeasDef1} is now $H:=L^{2}(%
\mathcal{R},\mathcal{H}^{h})$. The unitary operators $T$ and $U$ from (\ref%
{HausEq3}--\ref{HausEq4}) are now
\begin{equation}
(T_{k}f)(x):=f(x-k)\text{, \qquad }f\in H\text{, }x\in \mathcal{R}\text{, }%
k\in \mathbb{Z}^{d}  \label{IFSeq10}
\end{equation}%
and
\begin{equation}
(Uf)(x)=\frac{1}{\sqrt{\#(\mathcal{S})}}f(A^{-1}x)\text{, \qquad }f\in H%
\text{, }x\in \mathcal{R}\text{.}  \label{IFSeq11}
\end{equation}%
The commutation relation from proposition \ref{HausMeasProp3} in
its general form is
\begin{equation}
UT_{k}U^{-1}=T_{Ak}\text{,\qquad }k\in \mathbb{Z}^{d}\text{.}
\label{IFSeq12}
\end{equation}%
We now need the familiar duality between the two groups
$\mathbb{Z}^{d}$, and $\mathbb{T}^{d}=\mathbb{R}^{d}\diagup
\mathbb{Z}^{d}$, which identifies points $n\in \mathbb{Z}^{d}$
with monomials on $\mathbb{T}^{d}$ as follows,
\begin{equation}
z^{n}=z_{1}^{n_{1}}z_{2}^{n_{2}}\cdots z_{d}^{n_{d}}=e^{i2\pi
n_{1}\theta _{1}}e^{i2\pi n_{2}\theta _{2}}\cdots e^{i2\pi
n_{d}\theta _{d}}\text{.} \label{IFSeq13}
\end{equation}%
Note that (\ref{IFSeq13}) identifies the torus $\mathbb{T}^{d}$ with the $d$%
-cube
\begin{equation*}
\{(\theta _{1},\ldots ,\theta _{d})\mid 0\leq \theta _{i}<1\text{, }%
i=1,\ldots ,d\}\text{.}
\end{equation*}%
Since $\mathbf{C}$ is naturally identified with a subset of $\mathbb{T}^{d}$%
, we may view the monomials $\{z^{n}\mid n\in \mathbb{Z}^{d}\}$ as
functions on $\mathbf{C}$ by restriction. We say that the system
$(A,\mathcal{S})$ is
of \textit{orthogonal type} if there is a subset $\mathcal{T}$ of $\mathbb{Z}%
^{d}$ such that the set of functions $\{z^{n}\mid n\in
\mathcal{T}\}$ is an
orthonormal basis (ONB) in the Hilbert space $L^{2}(\mathbf{C},\mu _{(A,%
\mathcal{S})})$. If there is no subset $\mathcal{T}$ with this
ONB-property we say that $(A,\mathcal{S})$ is of
\textit{non-orthogonal type.} The authors of \cite{JoPe98} showed
that $(4,\{0,2\})$ is of orthogonal type, while $(3,\{0,2\})$ is
not. So for the Cantor set $\mathbf{C}_{4}$ there is
an ONB $\{z^{n}\mid n\in \mathcal{T}\}$ for a subset $\mathcal{T}$ of $%
\mathbb{Z}$; in fact we may take
\begin{equation}
\mathcal{T}=\{0,1,4,5,16,17,20,21,24,25,\cdots \}=\left\{ \dsum\limits_{0}^{%
\limfunc{finite}}n_{i}4^{i}\mid n_{i}\in \{0,1\}\right\} \text{.}
\label{IFSeq14}
\end{equation}%
For the middle-third Cantor set $\mathbf{C}_{3}$ it can be checked that $%
\{z^{n}\mid n\in \mathbb{Z\}}$ contains no more than two elements
which are orthogonal in $L^{2}(\mathbf{C}_{3},\mu _{3})$.

\begin{theorem}
\label{IFStheo1}Let $(A,\mathcal{S})$ be an affine IFS in
$\mathbb{R}^{d}$,
and suppose $\mathcal{S}$ has an extension to a set of $A$-residues in $%
\mathbb{Z}^{d}$. Let
\begin{equation*}
h=\frac{\log \#(\mathcal{S})}{\log \left\vert \det A\right\vert
}\text{,}
\end{equation*}%
and let $(\mathbf{C},\mu )$ be as above\emph{;} i.e., depending on $(A,%
\mathcal{S})$, and let $\mathcal{R}$ be defined from $\mathbf{C}$
in the usual way as in \textup{(\ref{IFSeq9})}. Assume further
that
\begin{equation}\label{eqIFSextra}
\mathbf{C}\cap(\mathbf{C}+k)=\emptyset,\quad(k\in\mathbb{Z}^d\setminus\{0\}).
\end{equation}
Then the system $(A,\mathcal{S})$
is of orthogonal type if and only iff there is a subset $\mathcal{T}$ in $%
\mathbb{Z}^{d}$ such that
$$
\left\{\left( \# (\mathcal{S}) \right)^{n/2} e^{i2\pi A^nk\cdot x}\chi _{\mathbf{C}}(A^nx-\ell )\mid k\in \mathcal{T}%
\text{, }(n=0\mbox{ and } \ell \in \mathbb{Z}^{d})\mbox{ or
}\right.$$
\begin{equation}
\left.(n\geq 1\mbox{ and }\ell\not\equiv s\mod A\mbox{ for all
}s\in\mathcal{S}) \right\}  \label{IFSeq15}
\end{equation}%
is an orthonormal basis in the Hilbert space $L^{2}(\mathcal{R},\mathcal{H}%
^{h})$.
\end{theorem}
\begin{remark}
The significance of the assumption (\ref{eqIFSextra}) is
illustrated in \cite{BrJo99}. Also note that (\ref{eqIFSextra}) is
automatically satisfied if $\mathbf{C}=\mathbf{C}(A,\mathcal{S})$
is contained in a $\mathbb{Z}^d$-tile. This is the case for the
example $(A,\mathcal{S})=(4,\{0,2\})$, but there are examples in
$d=2$ where it is not.
\end{remark}
\begin{proof}
A simple check shows that
$$\mathcal{R}=\bigcup\left\{A^{-n}(\mathbf{C}+l) \mid (n=0\mbox{ and } \ell \in \mathbb{Z}^{d})\right.
$$
$$\left.\mbox{ or }(n\geq 1\mbox{ and }\ell\not\equiv s\mod A\mbox{ for all }s\in\mathcal{S})\right\},$$
and the union is disjoint. Suppose $(A,\mathcal{S})$ is of
orthogonal type. We saw in Section 2 that the restriction of the
Hausdorff measure $\mathcal{H}^{h}$ to $\mathbf{C}$
agrees with the Hutchinson measure $\mu =\mu _{(A,\mathcal{S})}$ on $\mathbf{%
C=C}_{(A,\mathcal{S})}$. Hence density of $\{z^{n}\mid n\in
\mathcal{T}\}$
in $L^{2}(\mathbf{C},\mu )$ implies density of $\{e^{i2\pi k\cdot x}\chi _{%
\mathbf{C}}(x)\mid k\in \mathcal{T}\}$ in the subspace $L^{2}(\mathbf{C},%
\mathcal{H}^{h})$ of $L^{2}(\mathcal{R},\mathcal{H}^{h})$. Now the
formula for $\mathcal{R}$ implies that the functions in
(\ref{IFSeq15}) are dense in $L^2(\mathcal{R},\mathcal{H}^h)$.

Suppose conversely that the family (\ref{IFSeq15}) is dense in $L^{2}(%
\mathcal{R},\mathcal{H}^{h})$. Then $\{z^{n}\mid n\in
\mathcal{T}\}$ must be dense in $L^{2}(\mathbf{C},\mu )$ since
$\mathbf{C}$ is the support of Hutchinson's measure $\mu $, and
since $\mu $ restricts $\mathcal{H}^{h}$.
\end{proof}

\begin{corollary}
\label{IFScor1}Let $(\mathbf{C}_{4},\mu _{4})$ be the Cantor
construction in
the unit interval $I\cong \mathbb{T}^{1}$ defined by the IFS $\sigma _{0}(x)=%
\frac{x}{4},$ $\sigma _{2}(x)=\frac{x+2}{4}$\emph{;} i.e., by $(A,\mathcal{S}%
)=(4,\{0,2\})$, and let $\mathcal{R}$ be the subset of
$\mathbb{R}$ defined in \textup{(}\ref{IFSeq9}\textup{)}. Then the
family of functions
$$
\left\{ 2^{n/2}e^{i2\pi 4^nkx}\chi _{\mathbf{C}}(4^nx-\ell )\mid
k\in \{0,1,4,5,16,17,\cdots \},\right.$$
\begin{equation}\left.\ell \in
\left\{\begin{array}{ccc}
\mathbb{Z}&\mbox{if}&n=0\\
\mathbb{Z}\setminus(4\mathbb{Z}+\{0,2\})&\mbox{if}&n\geq1.
\end{array}\right.
\right\}  \label{IFSeq16}
\end{equation}%
forms an orthonormal basis in the Hilbert space $L^{2}(\mathcal{R},\mathcal{H%
}^{\frac{1}{2}})$.
\end{corollary}

\begin{proof}
This is a direct application of the theorem as the subset
\begin{equation*}
\mathcal{T}=\{0,1,4,5,16,17,\cdots \}
\end{equation*}%
from (\ref{IFSeq14}) and (\ref{IFSeq16}) satisfies the basis property for $%
\mathbf{C}_{4},\mu _{4}$ by Theorem 3.4 in \cite{JoPe98}.
\end{proof}

The next result makes clear the notion of gap-filling wavelets in
the context of iterated function systems (IFS). While it is stated
just for a particular example, the idea carries over to general
IFSs. Note that in the
system (\ref{IFSeq17}) below of wavelet functions, the two $\psi _{2}$ and $%
\psi _{3}$ are gap-filling.

\begin{corollary}
\label{IFScor2}Let $\mathbf{C}=\mathbf{C}_{4}$ be the Cantor set
determined from the IFS, $\sigma _{0}(x)=\frac{x}{4},$ $\sigma
_{2}(x)=\frac{x+2}{4}$, from the previous corollary. Then the
three functions
\begin{eqnarray}
\psi _{1}(x) &:&=\chi _{\mathbf{C}}(4x)-\chi _{\mathbf{C}}(4x-2)
\label{IFSeq17} \\
\psi _{2}(x) &:&=\sqrt{2}\chi _{\mathbf{C}}(4x-1)  \notag \\
\psi _{3}(x) &:&=\sqrt{2}\chi _{\mathbf{C}}(4x-3)  \notag
\end{eqnarray}%
generate an orthonormal wavelet basis in the Hilbert space $L^{2}(\mathcal{R}%
,\mathcal{H}^{\frac{1}{2}})$. Specifically, the family
\begin{equation}
\left\{ 2^{\frac{k}{2}}\psi _{i}(4^{k}x-\ell )\mid i=1,2,3\text{,
}k\in \mathbb{Z}\text{, }\ell \in \mathbb{Z}\right\}
\label{IFSeq18}
\end{equation}%
is an orthonormal basis in
$L^{2}(\mathcal{R},\mathcal{H}^{\frac{1}{2}})$.
\end{corollary}

\begin{proof}
We noted that our results in propositions \ref{HausMeasProp2} and
\ref{HausMeasProp4} apply more generally to IFSs of affine type.
So the result amounts to checking the general orthogonality
relations for the functions $m_{0},m_{1},m_{2},m_{3}$
on $\mathbb{T}$ which define wavelet filters for the system in (\ref{IFSeq18}%
). Note that from (\ref{IFSeq18}) the subband filters
$\{m_{i}\}_{i=0}^{3}$ are as follows, $z\in \mathbb{T}$:
\begin{eqnarray*}
m_{0}(z) &=&\frac{1}{\sqrt{2}}(1+z^{2}) \\
m_{1}(z) &=&\frac{1}{\sqrt{2}}(1-z^{2}) \\
m_{2}(z) &=&z \\
m_{3}(z) &=&z^{3}\text{.}
\end{eqnarray*}%
Since the $4\times 4$ matrix in the system
\begin{equation*}
\left(
\begin{array}{c}
m_{0}(z) \\
m_{1}(z) \\
m_{2}(z) \\
m_{3}(z)%
\end{array}%
\right) =\left(
\begin{array}{cccc}
\frac{1}{\sqrt{2}} & 0 & \frac{1}{\sqrt{2}} & 0 \\
\frac{1}{\sqrt{2}} & 0 & -\frac{1}{\sqrt{2}} & 0 \\
0 & 1 & 0 & 0 \\
0 & 0 & 0 & 1%
\end{array}%
\right) \left(
\begin{array}{c}
1 \\
z \\
z^{2} \\
z^{3}%
\end{array}%
\right)
\end{equation*}%
is clearly unitary, the result follows from a direct computation;
see also the proof of theorem \ref{XformTheo1}.

To verify that the Ruelle operator $R=R_{m_0}$ given by
\begin{eqnarray*}
(Rf)(z)
&=&\frac{1}{4}\dsum\limits_{w^{4}=z}\left\vert
m_{0}(w)\right\vert ^{2}f(w) \\
&=&\frac{1}{4}\dsum\limits_{w^{4}=z}\left(
1+\frac{w^{2}+w^{-2}}{2}\right) f(w)
\end{eqnarray*}%
satisfies the two conditions

(a) $\dim \{f\in C(\mathbb{T})\mid Rf=f\}=1$, and

(b) for all $\lambda \in \mathbb{C}$, $\left\vert \lambda
\right\vert =1$, and $\lambda \neq 1$, $\dim \{f\in
C(\mathbb{T})\mid Rf=\lambda f\}=0$, we may again apply the
theorem from \cite{Nus98} or the results of section 6 below.
\end{proof}

For the more general affine IFSs the results above extend as
follows.

Consider the affine IFS $(\sigma _{i})_{i=1}^{p}$ with
\begin{equation*}
\sigma _{i}(x)=\frac{1}{N}(x+a_{i}),\quad (x\in \mathbb{R}),
\end{equation*}%
where $N\geq 2$ is an integer and $(a_{i})_{i=1}^{p}$ are distinct
integers in $\{0,...,N-1\}$. Then by \cite{Fal}, there is a unique
compact subset $K$ of $\mathbb{R}$ which is the attractor of the
IFS, i.e.,
\begin{equation*}
\mathbf{C}=\cup _{i=1}^{p}\sigma _{i}(\mathbf{C}).
\end{equation*}%
Actually, one can give a more explicit description of this
attractor, namely
\begin{equation*}
\mathbf{C}=\{\sum_{j\geq 1}d_{j}N^{-j}\,|\,d_{j}\in \{a_{1},...,a_{p}\}\text{%
, }j\geq 1\}.
\end{equation*}%
Since the digits $a_{i}$ are distinct and less then $N$, $K$ is
contained in $[0,1]$, and the sets $\sigma_{i}(K)$ are almost
disjoint (they have at most one point in common, those of the form
$k/N$ for some $k\in \{1,...,N-1\}$.

The Hausdorff dimension of $K$ is $\log_Np$.

Now consider the set
\begin{equation*}
\mathcal{R}=\{\sum_{j\geq -m}d_{j}N^{-j}\,|\,m\in
\mathbb{Z},d_{j}\in \{a_{1},...,a_{p}\}\text{ for all but finitely
many indices }j\}
\end{equation*}%
$\mathcal{R}$ is invariant under integer translations
\begin{equation*}
\mathcal{R}+k=\mathcal{R},\quad (k\in \mathbb{Z}),
\end{equation*}%
and it is invariant under dilation by $N$
\begin{equation*}
N\mathcal{R}=\mathcal{R}.
\end{equation*}

Endow $\mathcal{R}$ with the Hausdorff measure $\mathcal{H}^s$ for $%
s=\log_Np $, and on $L^2(\mathcal{R},\mathcal{H}^s)$, define the
translation operator
\begin{equation*}
Tf(x)=f(x-1),\quad(x\in\mathcal{R},f\in
L^2(\mathcal{R},\mathcal{H}^s)),
\end{equation*}
and the dilation operator
\begin{equation*}
Uf(x)=\sqrt{\frac{1}{p}}f\left(\frac{x}{N}\right),\quad(x\in\mathcal{R},f\in
L^2(\mathcal{R},\mathcal{H}^s)).
\end{equation*}
These are unitary operators satisfying the commutation relation
\begin{equation*}
UTU^{-1}=T^N.
\end{equation*}

Let $\varphi :=\chi _{\mathbf{C}}$. The function $\varphi $ is an
orthogonal scaling function for
$L^{2}(\mathcal{R},\mathcal{H}^{s})$, with filter
\begin{equation}\label{eq3_19_1}
m_{0}(z)=\sqrt{\frac{1}{p}}\sum_{i=1}^{p}z^{a_{i}},
\end{equation}%
so it satisfies the following conditions:

\begin{enumerate}
\item \textbf{[Orthogonality]}
\begin{equation*}
\left\langle T^{k}\varphi \mid \varphi \right\rangle =\delta
_{k},\quad (k\in \mathbb{Z}).
\end{equation*}

\item \textbf{[Scaling equation]}
\begin{equation*}
U\varphi=\sum_{i=1}^p\sqrt{\frac{1}{p}}T^{a_i}\varphi=m_0(T).
\end{equation*}

\item \textbf{[Cyclicity]}
\begin{equation*}
\overline{\func{span}}\{U^{-n}T^{k}\varphi \,|\,n,k\in \mathbb{Z}\}=L^{2}(%
\mathcal{R},\mathcal{H}^{s}).
\end{equation*}
\end{enumerate}

Next, we define the wavelets. For this, we need the "high-pass" filters $%
m_1,..., m_{N-1}$ such that the matrix
\begin{equation*}
\frac{1}{\sqrt{N}}(m_i(\rho^jz))_{i,j=0}^{N-1},
\end{equation*}
is unitary for almost every $z$. ($\rho=e^{2\pi i/N}$).

First, we define the filters for the gap-filling wavelets $%
\psi_1,...,\psi_{N-p}$. The set
$G=\{0,...,N-1\}\setminus\{a_1,...,a_p\}$
has $N-p$ elements. We label the functions $z\mapsto z^d$ for $d\in G$, by $%
m_1,...,m_{N-p}$.

The remaining $p-1$ filters are for the detail-filling wavelets. Let $%
\eta=e^{2\pi i/p}$. Define
\begin{equation*}
m_{N-p+k}(z)=\sqrt{\frac{1}{p}}\sum_{i=1}^p\eta^{k(i-1)}z^{a_i},\quad(k\in%
\{1,...,p-1\}).
\end{equation*}

We have to check that
\begin{equation}  \label{equnit1}
\frac{1}{N}\sum_{w^N=z}m_i(w)\overline{m}_j(w)=\delta_{ij},\quad(z\in\mathbb{%
T},i,j\in\{0,...,N-1\}).
\end{equation}
For this we use the following identity:
\begin{equation*}
\sum_{w^N=z}w^k=0,\quad(z\in\mathbb{T},k\not\equiv 0\mod N).
\end{equation*}
Therefore, if $f_1(z)=\sum_{i=0}^{N-1}\alpha_iz^i$, $f_2=\sum_{i=0}^{N-1}%
\beta_iz^i$, then
\begin{equation*}
\frac{1}{N}\sum_{w^N=z}f_1(w)\overline{f_2(w)}=\frac{1}{N}%
\sum_{i,j=0}^{N-1}\alpha_i\overline{\beta}_j\sum_{w^N=z}w^{i-j}=%
\sum_{i=0}^{N-1}\alpha_i\overline{\beta}_j.
\end{equation*}
Applying these to the filters $m_i, (i\in\{0,...,N-1\})$, we obtain (\ref%
{equnit1}).

With these filters, we construct the wavelets in the usual way:
\begin{equation*}
\psi _{i}=U^{-1}m_{i}(T)\varphi ,\quad (i\in \{1,...,N-1\}),
\end{equation*}%
and
\begin{equation*}
\{U^{m}T^{n}\psi _{i}\,|\,m,n\in \mathbb{Z},i\in \{1,...,N-1\}\}
\end{equation*}%
is an orthonormal basis for $L^{2}(\mathcal{R},\mathcal{H}^{s}).$

Let $N\in \mathbb{Z}_{+}$ be as above, and consider $\mathcal{S}%
=\{a_{1},\cdots ,a_{p}\}\subset \{0,1,2,\cdots ,N-1\}$. A second subset $%
\mathcal{B}=\{b_{1},\cdots ,b_{p}\}\subset \mathbb{Z}$ is an
$N$-dual if the $p\times p$ matrix
\begin{equation}
M_{N}(\mathcal{S},\mathcal{B})=\frac{1}{\sqrt{p}}\left( \exp
\left( i\frac{2\pi a_{j}b_{k}}{N}\right) \right) _{1\leq j,k\leq
p} \label{IFSeq20}
\end{equation}%
is unitary. When $N$ and $\mathcal{S}$ are given as specified, it
is not always true
that there is a subset $\mathcal{B}\subset \mathbb{Z}$ for which $M_{N}(%
\mathcal{S},\mathcal{B})$ is unitary. If for example $N=3$ and $\mathcal{S}%
=\{0,2\}$, then no $\mathcal{B}$ exists, while for $N=4$ and $\mathcal{S}%
=\{0,2\}$, we may take $\mathcal{B}=\{0,1\}$, and
\begin{equation*}
M_{4}(\mathcal{S},\mathcal{B})=\frac{1}{\sqrt{2}}\left(
\begin{array}{cc}
1 & 1 \\
1 & -1%
\end{array}%
\right)
\end{equation*}%
is of course unitary.

\begin{lemma}
\label{IFSlem1}\cite{JoPe98} Let $N$ and $\mathcal{S}$ be as
specified above, and suppose
\begin{equation*}
\mathcal{B}=\{b_{1},\cdots ,b_{p}\}\subset \mathbb{Z}
\end{equation*}
is an $N$-dual subset. Suppose $0\in \mathcal{B}$, and set
\begin{equation}
\Lambda =\Lambda _{N}(\mathcal{B}):=\left\{ \sum_{i=0}^{\limfunc{finite}%
}n_{i}N^{i}\mid n_{i}\in \mathcal{B}\right\} \text{.}
\label{IFSeq21}
\end{equation}%
Let $(\mathbf{C},\mu )=(\mathbf{C}_{(N,\mathcal{S})},\mu
_{(N,\mathcal{S})})$ be the Hutchinson pair. Then the set of
functions $\{z^{n}\mid n\in \Lambda \}$ is orthogonal in
$L^{2}(\mathbf{C},\mu )$\emph{;} i.e.,
\begin{equation}
\dint\limits_{\mathbf{C}}z^{n-n^{\prime }}d\mu (z)=\delta
_{n,n^{^{\prime }}},\qquad n,n^{\prime }\in \Lambda
\label{IFSeq22}
\end{equation}%
where we identify $\mathbf{C}$ as a subset of $\mathbb{T}^{1}$ via
\begin{equation*}
\mathbf{C}\ni \theta \longrightarrow e^{i2\pi \theta }\in \mathbb{T}^{1}%
\text{.}
\end{equation*}
\end{lemma}

\begin{proof}
Set $e(\theta )=e^{i2\pi \theta }$, and for $k\in \mathbb{R}$
\begin{equation}
B(k):=\dint\limits_{\mathbf{C}}e(k\theta )d\mu (\theta )\text{.}
\label{IFSeq23}
\end{equation}%
Using (\ref{IFSeq4}), we get
\begin{equation}
B(k)=\frac{1}{\sqrt{p}}m_{0}\left( \frac{k}{N}\right) B\left( \frac{k}{N}%
\right),  \label{IFSeq24}
\end{equation}%
where $m_0$ is defined in (\ref{eq3_19_1}).\par If $n,n^{\prime
}\in \Lambda $, and $n\neq n^{\prime }$, we get the representation
\begin{equation*}
n^{\prime }-n=b^{\prime }-b+mN^{\ell }\text{,\qquad }b,b^{\prime
}\in \mathcal{B}\text{, }m,\ell \in \mathbb{Z}\text{, }\ell \geq
1\text{.}
\end{equation*}%
As a result, the inner product in $L^{2}(\mathbf{C},\mu )$ is
\begin{equation}
\left\langle z^{n}\mid z^{n^{\prime }}\right\rangle _{\mu }=B(n^{\prime }-n)=%
\frac{1}{\sqrt{p}}m_{0}\left( \frac{b^{\prime }-b}{N}\right) B\left( \frac{%
n^{\prime }-n}{N}\right) \text{.}  \label{IFSeq25}
\end{equation}%
Since the matrix $M_{N}(\mathcal{S},\mathcal{B})$ is unitary,
\begin{equation*}
m_{0}\left( \frac{b^{\prime }-b}{N}\right) =0
\end{equation*}%
when $b^{\prime }\neq b$ in $\mathcal{B}$, and the result follows.
\end{proof}

Even if the matrix $M_{N}(\mathcal{S},\mathcal{B})$ is unitary,
the orthogonal functions $\{z^{n}\mid n\in \Lambda \}$ might not
form a basis for $L^{2}(\mathbf{C},\mu )$. From \cite{JoPe98}, we
know that it is an orthonormal basis (ONB) if and only if
\begin{equation}
\dsum\limits_{n\in \Lambda }\left\vert B(\xi -n)\right\vert ^{2}=1\text{%
\qquad a.e. }\xi \in \mathbb{R}\text{.}  \label{IFSeq26}
\end{equation}

Introducing the function
\begin{equation}
\Omega (\xi ):=\frac{1}{p}\dsum\limits_{b\in
\mathcal{B}}\left\vert B(\xi -n)\right\vert ^{2}\text{,}
\label{IFSeq27}
\end{equation}%
and the dual Ruelle operator
\begin{equation}
(R_{\mathcal{B}}f)(\xi ):=\frac{1}{p}\dsum\limits_{b\in \mathcal{B}%
}\left\vert m_{0}\left( \frac{\xi -b}{N}\right) \right\vert
^{2}f\left( \frac{\xi -b}{N}\right) \text{,}  \label{IFSeq28}
\end{equation}%
we easily verify that $\Omega $ and the constant function
$\hat{1}$ both solve the eigenvalue problem
$R_{\mathcal{B}}(f)=f$, both functions $\Omega $ and $\hat{1}$ are
continuous on $\mathbb{R}$, even analytic.

\begin{theorem}
\label{IFStheo2}If the space
\begin{equation}
\left\{ f\in \limfunc{Lip}(\mathbb{R)}\mid f\geq 0\text{, }f(0)=1\text{%
,\quad }R_{\mathcal{B}}(f)=f\right\}  \label{IFSeq29}
\end{equation}%
is one-dimensional, then $\Lambda (=\Lambda _{N}(\mathcal{B}))$
induces an
ONB\emph{;} i.e., $\{z^{n}\mid n\in \Lambda \}$ is an ONB in $L^{2}(\mathbf{C%
},\mu )$.
\end{theorem}

\begin{proof}
The result follows from the discussion and the added observation that $%
\Omega (0)=1$. This normalization holds since $0\in \mathcal{B}$
was assumed, and so $\left\langle e_{0}\mid e_{n}\right\rangle
_{\mu }=0$ for all $n\in \Lambda \diagdown \{0\}$.
\end{proof}

\begin{definition}
\label{IFSdef1}A $\mathcal{B}$-cycle is a finite set
$\{z_{1},z_{2},\ldots ,z_{k+1}\}\subset \mathbb{T}$, with a
pairing of points in $\mathcal{B}$, say $b_{1},b_{2},\ldots
,b_{k+1}\in \mathcal{B}$, such that
\begin{equation}
z_{i}=\sigma _{-b_{i}}(z_{i+1})\text{,\qquad
}z_{k+1}=z_{1}\text{,} \label{IFSeq30}
\end{equation}%
and $\left\vert m_{0}(z_{i})\right\vert ^{2}=p$. Equivalently, a
$\mathcal{B} $-cycle may be given by $\{\xi _{1},\ldots ,\xi
_{k+1}\}\subset \mathbb{R}$ satisfying
\begin{eqnarray*}
\xi _{i+1} &\equiv &b_{i}+N\xi _{i}\;\limfunc{mod}N\mathbb{Z} \\
\left( N^{k}-1\right) \xi _{1} &\equiv &b_{k}+Nb_{k-1}+\cdots +N^{k-1}b_{1}\;%
\limfunc{mod}N^{k}\mathbb{Z}\text{.}
\end{eqnarray*}
\end{definition}

\begin{theorem}
\label{IFStheo3}Let $N\in \mathbb{Z}_{+}$, $N\geq 2$ be given. Let $\mathcal{%
S}\subset \{0,1,\cdots ,N-1\}$, and suppose there is a
$\mathcal{B}\subset
\mathbb{Z}$ such that $0\in \mathcal{B}$, $\#(\mathcal{S})=\#(\mathcal{B})=p$%
, and the matrix
\begin{equation*}
M_{N}(\mathcal{S},\mathcal{B})=\frac{1}{\sqrt{p}}\left( \exp \left( i\frac{%
2\pi ab}{N}\right) \right)
\end{equation*}%
is unitary. Then $\{z^{n}\mid n\in \Lambda _{N}(\mathcal{B})\}$ is
an ONB for $L^{2}(\mathbf{C},\mu )$ where $\Lambda
_{N}(\mathcal{B})$ is defined in \textup{(}\ref{IFSeq21}\textup{)}
if the only $\mathcal{B}$-cycles are the singleton $\{1\}\subset
\mathbb{T}$.
\end{theorem}

\begin{proof}
By Theorem \ref{IFStheo2}, we need only verify that the absence of $\mathcal{%
B}$-cycles of order $\geq 2$ implies that the Perron-Frobenius eigenspace (%
\ref{IFSeq28}) is one-dimensional. But this follows from \cite[ Theorem 5.5.4%
]{BrJo02}. In fact, the argument from Chapter 5 in \cite{BrJo02}
shows that
the absence of $\mathcal{B}$-cycles of order $\geq 2$ implies that the $%
\mathcal{B}$-Ruelle operator $R_{\mathcal{B}}$ with $\sigma _{-b}(\xi ):=%
\frac{\xi -b}{N}$,
\begin{equation*}
(R_{\mathcal{B}}f)(\xi )=\frac{1}{p}\dsum\limits_{b\in
\mathcal{B}}\left( \left\vert m_{0}(\sigma _{-b}(\xi ))\right\vert
^{2}f(\sigma _{-b}(\xi ))\right)
\end{equation*}%
satisfies the two Perron-Frobenius properties:

\noindent (i) the only bounded continuous solutions $f$ to $R_{\mathcal{B}%
}(f)=f$ are the multiples of $\hat{1}$, and

\noindent (ii) for all $\lambda \in \mathbb{T}\diagdown \{1\}$,
the eigenvalue problem $R_{\mathcal{B}}(f)=\lambda f$ has no
non-zero bounded continuous solutions.
\end{proof}

\begin{example}
\label{IFSex1}\textup{(}An Application\textup{)} Let $N=4$, $\mathcal{S}%
=\{0,2\}$, and $\mathcal{B}=\{0,1\}$. Then
\begin{eqnarray*}
M_{4}(\mathcal{S},\mathcal{B}) &=&\frac{1}{\sqrt{2}}\left(
\begin{array}{cc}
1 & 1 \\
1 & -1%
\end{array}%
\right) \text{,} \\
\Lambda _{4}(\mathcal{B}) &=&\{0,1,4,5,16,17,20,21,\cdots \}\text{, and} \\
(R_{\mathcal{B}}f)(\xi ) &=&\cos ^{2}(2\pi \xi )\;f\left( \frac{\xi }{4}%
\right) +\sin ^{2}(2\pi \xi )\;f\left( \frac{\xi -1}{4}\right)
\end{eqnarray*}%
and there is only on $\mathcal{B}$-cycle, the singleton
$\{1\}\subset
\mathbb{T}$. Recall from \cite{Hut81} that the Hutchinson construction of $(\mathbf{C},\mu )$ identifies $%
\mathbf{C}$ as the Cantor set arising by the subdividing algorithm
starting with the unit interval $I$ dividing into four equal
subintervals and dropping the second and the fourth at each step
in the algorithm. The measure $\mu $ is the restriction of
$\mathcal{H}^{\frac{1}{2}}$ to $\mathbf{C}$, and it
follows from the last theorem that $\{z^{n}\mid n\in \Lambda _{4}(\mathcal{B}%
)\}$ is an ONB for $L^{2}(\mathbf{C},\mu )$. The dual system
$\{\sigma
_{-b}\mid b\in \mathcal{B}\}$; i.e., $\sigma _{0}(\xi )=\frac{\xi }{4}$, $%
\sigma _{-1}(\xi )=\frac{\xi -1}{4}$, generates a Cantor subset $\mathbf{C}_{%
\mathcal{B}}\subset \lbrack -1,0]$ also of Hausdorff dimension $\frac{1}{2}$%
. Note that the fractional version of the Ruelle operator
$R_{\mathcal{B}}$ does not map $1$-periodic functions into
themselves; in general
\begin{equation*}
(R_{\mathcal{B}}f)(\xi )\neq (R_{\mathcal{B}}f)(\xi +1)\text{;}
\end{equation*}%
in fact
\begin{equation*}
(R_{\mathcal{B}}f)(\xi +1)=\cos ^{2}(2\pi \xi )f\left( \frac{\xi +1}{4}%
\right) +\sin ^{2}(2\pi \xi )f\left( \frac{\xi }{4}\right)
\text{;}
\end{equation*}%
so
\begin{equation*}
R_{\mathcal{B}}f(\xi )=R_{\mathcal{B}}f(\xi +1)
\end{equation*}%
holds only if
\begin{equation*}
\cos ^{4}(2\pi \xi )f\left( \frac{\xi +1}{4}\right) =\sin
^{4}(2\pi \xi )f\left( \frac{\xi -1}{4}\right) \text{.}
\end{equation*}
\end{example}

The following tables of similar examples is included hopefully
offering the reader a glimpse of the variety of examples, all of
orthogonal type. The tables also offers some insight into the
duality between the two systems,
one in the $x$-variable and the other in the Fourier dual variable $\xi $%
\begin{equation*}
\begin{array}{cccccc}
N & p & \mathcal{S} & \mathcal{B} & M_{N}(\mathcal{S},\mathcal{B}) & \text{%
Hausdorff Dim}\emph{.} \\
4 & 2 & \{0,2\} & \{0,1\} & \frac{1}{\sqrt{2}}\left(
\begin{array}{ll}
1 & 1 \\
1 & -1%
\end{array}%
\right) _{\mathstrut } & \frac{1}{2} \\
6 & 2 & \{0,3\} & \{0,1\} & \frac{1}{\sqrt{2}}\left(
\begin{array}{ll}
1 & 1 \\
1 & -1%
\end{array}%
\right) _{\mathstrut } & \log _{6}(2) \\
6 & 2 & \{0,1\} & \{0,3\} & \frac{1}{\sqrt{2}}\left(
\begin{array}{ll}
1 & 1 \\
1 & -1%
\end{array}%
\right) _{\mathstrut } & \log _{6}\left( 2\right) \\
6 & 3 & \{0,2,4\} & \{0,1,2\} &
\begin{array}{l}
\frac{1}{\sqrt{3}}\left(
\begin{array}{lll}
1 & 1 & 1 \\
1 & \zeta _{3} & \zeta _{3}^{2} \\
1 & \zeta _{3}^{2} & \zeta _{3}%
\end{array}%
\right) \\
\text{where }\zeta _{3}=\exp \left( i\frac{2\pi }{3}\right)%
\end{array}
& \log _{6}\left( 3\right)%
\end{array}%
\end{equation*}

\begin{equation*}
\begin{array}{ccc}
N & p & \Lambda _{N}(\mathcal{B}) \\
4 & 2 & \{0,1,4,5,16,17,20,21,\cdots \} \\
6 & 2 & \{0,1,6,7,36,37,42,43,\cdots \} \\
6 & 2 & \{0,3,6,9,36,39,42,45,\cdots \} \\
6 & 3 & \{0,1,2,6,7,8,36,37,38,42,43,44,\cdots \}%
\end{array}%
\end{equation*}

\begin{equation*}
\begin{array}{ccc}
N & p & (R_{B}f)(\xi ) \\
4 & 2 & \cos ^{2}(2\pi \xi )f\QOVERD( ) {\xi }{4}+\sin ^{2}(2\pi
\xi
)f\QOVERD( ) {\xi -1}{4}_{\mathstrut } \\
6 & 2 & \cos ^{2}\QOVERD( ) {\pi \xi }{2}f\QOVERD( ) {\xi
}{6}+\sin
^{2}\QOVERD( ) {\pi \xi }{2}f\QOVERD( ) {\xi -1}{6}_{\mathstrut } \\
6 & 2 & \cos ^{2}\QOVERD( ) {\pi \xi }{6}f\QOVERD( ) {\xi
}{6}+\sin
^{2}\QOVERD( ) {\pi \xi }{6}f\QOVERD( ) {\xi -3}{6}_{\mathstrut } \\
6 & 3 & W\QOVERD( ) {\xi }{6}f\QOVERD( ) {\xi }{6}+W\QOVERD( )
{\xi -1}{6}f\QOVERD( ) {\xi -1}{6}+W\QOVERD( ) {\xi
-2}{6}f\QOVERD( ) {\xi
-2}{6}_{\mathstrut } \\
&  & \text{where }W(\xi ):=\frac{2\cos ^{2}(2\pi \xi )-\sin ^{2}(3\pi \xi )}{%
3}\text{.}%
\end{array}%
\end{equation*}

\section{\label{ZakTransform}A Generalized Zak-Transform}

The notion of filter is imported into math from signal processing. It has
now been well adapted to wavelet analysis: In the familiar dyadic case the
two wavelet functions $\varphi $ (the father function) and $\psi $ (the
mother function) are in the Hilbert space $L^{2}(\mathbb{R)}$. In the $N$%
-adic case, the wavelet functions are $\varphi ,\psi _{1},\cdots ,\psi
_{N-1} $, and there are known conditions for when these functions are in $%
L^{2}(\mathbb{R)}$. The starting point is the scaling identity (\ref%
{IntroEq3c}) satisfied by $\varphi $. Introducing the wavelet filter $%
m_{0}(z)=\dsum\limits_{k}a_{k}z^{k}$ as a function on $\mathbb{T=R}/2\pi
\mathbb{Z}$, and the transfer operator $R_{m_{0}}$ in (\ref{IntroEq3g}), we
note that necessary conditions for $\varphi $ to be in $L^{2}(\mathbb{R)}$
are $\left\vert m_{0}(1)\right\vert =\sqrt{N}$, the low-pass condition, and $%
R_{m_{0}}(\hat{1})=\hat{1}$, where $\hat{1}$ is the constant
function $1$. Let $\delta _{1}$ denote the Dirac mass at $z=1$. It
follows that $\delta _{1}R_{m_{0}}=\delta _{1}$. But what if, for
some $m_0$, $\delta _{1}$ does not satisfy this, so called
low-pass condition; but rather there is
some other probability measure $\nu $ on $\mathbb{T}$ which is $R_{m_{0}}$%
-invariant; i.e., satisfies $\nu R_{m_{0}}=\nu $, and which is
singular with full support. A main point in our paper is that this
alternative introduces fractal analysis into the wavelet
construction.

Traditionally, the Zak-transform \cite{Dau92} is a standard tool of analysis
in $L^{2}(\mathbb{R})$, and in this section it is extended to the abstract
case of Hilbert space. If $\{T_{k}$:~$k\in \mathbb{Z}\}$ denotes translation
$(T_{k}f)(x)=f(x-k)$, $x\in \mathbb{R}$, $k\in \mathbb{Z}$, $f\in L^{2}(%
\mathbb{R})$, we set%
\begin{equation*}
(Zf)(z,x)=\sum_{k\in \mathbb{Z}}z^{k}T_{k}f(x),\qquad z\in
\mathbb{T}\text{,} x\in I,
\end{equation*}%
and we check that $Z$ defines a unitary isomorphism of $L^{2}(\mathbb{R})$
onto $L^{2}(\mathbb{T}\times I)$ where $\mathbb{T}$ is the torus, and $I$
the unit-interval $I=[0,1)$. The measure on $\mathbb{T}$ is Haar measure,
denoted $\mu $; i.e.,
\begin{equation*}
\dint\limits_{\mathbb{T}}\cdots d\mu (z)=\frac{1}{2\pi }\dint\limits_{0}^{2%
\pi }\cdots d\theta \text{,\qquad where }z=e^{i\theta }\text{.}
\end{equation*}

Let $\mathcal{H}$ be a Hilbert space, $U$ a unitary operator in $\mathcal{H}$%
, and $T$:~$\mathbb{Z}\rightarrow \mathcal{U}(\mathcal{H)}$ a
unitary representation. Let $N\geq 2$, and suppose that
\begin{equation}
UT_{k}U^{-1}=T_{Nk},\qquad k\in \mathbb{Z}\text{.}  \label{ZakEq0}
\end{equation}

In general, the inner product in a Hilbert space $\mathcal{H}$ will be
written $\left\langle \cdot \mid \cdot \right\rangle $. If $f\in \mathcal{H}$%
, then the form $g\rightarrow \left\langle f\mid g\right\rangle $ is taken
linear on $\mathcal{H}$; and $\left\langle f\mid f\right\rangle =\left\Vert
f\right\Vert ^{2}$.

For $f_{i}\in \mathcal{H},i=1,2$, we introduce the following function%
\begin{equation}
p(f_{1,}f_{2})(z)=\sum_{k\in \mathbb{Z}}z^{k}\left\langle T_{k}f_{1}\mid
f_{2}\right\rangle  \label{ZakEq1}
\end{equation}%
defined formally for $z\in \mathbb{T}$. Let $m_{0}\in L^{\infty }(\mathbb{T}%
) $ be given, and suppose that
\begin{equation}
\frac{1}{N}\dsum\limits_{\underset{w^{N}=z}{{\small w\in }\mathbb{T}}%
}\left\vert m_{0}(w)\right\vert ^{2}=1,\qquad \text{a.e.}\qquad z\in \mathbb{%
T}\text{.}  \label{ZakEq2}
\end{equation}

Note that the sum in (\ref{ZakEq2}) is finite, since for each $z\in \mathbb{T%
}$, the equation $w^{N}=z$ has precisely $N$ solutions. In fact, the cyclic
group $\mathbb{Z}_{N}$ acts transitively on this set of solutions $\{w\}$.

(We will work with the torus $\mathbb{T}$ in anyone of its three familiar
incarnations: (i) $\left\{ z\in \mathbb{C}\mid \left\vert z\right\vert
=1\right\} $, (ii) the quotient group $\mathbb{R}/2\pi \mathbb{Z}$, or (iii)
the period interval $[0,2\pi )$ via the identification $z=e^{i\theta }.$
With this identification, we get the familiar description of the set $%
\left\{ w\in \mathbb{T}\mid w^{N}=z(=e^{i\theta })\right\} $ in the
summation on the left-hand side of (\ref{ZakEq2}) as the $N$ distinct
frequency bands $\left\{ \frac{\theta +2\pi k}{N}\mid k=0,1,\cdots
,N-1\right\} $. The sub-interval $[0,\frac{2\pi }{N})$ represents the low
frequency band.)

The operator $m_{0}(T)$ is defined from the spectral theorem in the usual
way: If the spectral measure of $T$ is denoted $E_{T}$, then $E_{T}$ is a
projection valued measure on $\mathbb{T}$, and we have the following three
identities:%
\begin{equation}
\left\Vert f\right\Vert ^{2}=\dint\limits_{\mathbb{T}}\left\Vert
E_{T}(dz)f\right\Vert ^{2},\qquad f\in \mathcal{H}\text{,}  \label{ZakEq3}
\end{equation}

\begin{equation}
T_{k}=\dint\limits_{\mathbb{T}}z^{k}E_{T}(dz),\qquad k\in Z\text{,}
\label{ZakEq4}
\end{equation}%
and by functional calculus,%
\begin{equation}
m_{0}(T)=\dint\limits_{\mathbb{T}}m_{0}(z)E_{T}(dz)  \label{ZakEq5}
\end{equation}%
If%
\begin{equation}
m_{0}(z)=\dsum\limits_{k\in \mathbb{Z}}a_{k}z^{k}  \label{ZakEq6}
\end{equation}%
is the Fourier series of $m_{0}$, it follows that $m_{0}(T)=\sum_{k\in
\mathbb{Z}}a_{k}T_{k}$ is then well defined.

The following operator $R=R_{m_{0}}$, called the \emph{Ruelle operator},
(see \cite{Jorgen01}) is acting on functions $h$ or $\mathbb{T}$ as follows,%
\begin{equation}
(Rh)(z)=\frac{1}{N}\dsum\limits_{\underset{w^{N}=z}{{\small w\in }\mathbb{T}}%
}\left\vert m_{0}(w)\right\vert ^{2}h(w),\qquad z\in \mathbb{T}\text{.}
\label{ZakEq7}
\end{equation}%
Let $\hat{1}$ denote the constant function $1$ on $\mathbb{T}$. Then
condition (\ref{ZakEq2}) amounts to the eigenvalue equation
\begin{equation}
R(\hat{1})=\hat{1}  \label{ZakEq8}
\end{equation}

On the Hilbert space $\mathcal{H}$, we introduce the operator%
\begin{equation}
M:=U^{-1}m_{0}(T)\text{.}  \label{ZakEq9}
\end{equation}%
It is called the \emph{cascade approximation operator}. In the special case
when $T_{k}f(x)=f(x-k)$, $(Uf)(x)=N^{-\frac{1}{2}}f(\frac{x}{N})$, and $%
\mathcal{H}=L^{2}(\mathbb{R})$, then%
\begin{equation}
(Mf)(x)=\sqrt{N}\dsum\limits_{k}a_{k}f(Nx-k)  \label{ZakEq10}
\end{equation}%
where $\{a_{k}$:~$k\in \mathbb{Z}\}$ is the sequence of Fourier coefficients
of $m_{0}$, see (\ref{ZakEq6}). In this case $M$ is also called \textit{the}
\textit{wavelet subdivision operator}. The following general lemma applies
also to the case of wavelets in $L^{2}(\mathbb{R})$. The advantage of (\ref%
{ZakEq9}) over (\ref{ZakEq10}) is that (\ref{ZakEq9}) is defined for all
systems $U,T$ satisfying (\ref{ZakEq0}) and applies in particular to our
present fractal examples.

For measurable functions $\xi $ and $\eta $ on $\mathbb{T}$, formula (\ref%
{ZakEq5}) represents the usual functional calculus; i.e., $\xi (T)=\int_{%
\mathbb{T}}\xi (z)E_{T}(dz)$. Setting $\pi (\xi ):=\xi (T)$, we get $\pi
(\xi \eta )=\pi (\xi )\pi (\eta )$, $\pi (\hat{1})=I=$ the identity
operator, $\pi (\,\overline{\xi }\,)=\pi (\xi )^{\ast }=$ the adjoint
operator. These properties together state that $\pi (=\pi _{E_{T}})$ defines
a $\ast $-representation of $L^{\infty }(\mathbb{T)}$ acting on the Hilbert
space $\mathcal{H}_{T}$ of the translation operators $\{T_{k}$:~$k\in
\mathbb{Z}\}$.

\begin{lemma}
\label{ZakLemma1}Let $\mathcal{H}$, $T$, $U$, and $m_{0}$ be as described
above, and let the operators $M$ and $R$ be the corresponding operators;
i.e., the cascade operator, and Ruelle operator, respectively. Then the
identity%
\begin{equation}
R(p(f_{1},f_{2}))=p(Mf_{1},Mf_{2})  \label{ZakEq11}
\end{equation}%
holds for all $f_{1,}f_{2}\in \mathcal{H}$, where the two sides in \textup{(}%
\ref{ZakEq11}\textup{)} are viewed as functions on $\mathbb{T}$.
\end{lemma}

\begin{proof}
Let $\xi \in C(\mathbb{T})$. Then it follows from (\ref{ZakEq1}) and (\ref%
{ZakEq6}) that
\begin{equation}
\dint\limits_{\mathbb{T}}\xi (z)p(f_{1},f_{2})(z)d\mu (z)=\left\langle
f_{1}\mid \xi (T)f_{2}\right\rangle  \label{ZakEq12}
\end{equation}%
where $\mu $ is the Haar measure on $\mathbb{T}$, $\left\langle \cdot \mid
\cdot \right\rangle $ is the inner product of $\mathcal{H}$, and $\xi
(T)=\int_{\mathbb{T}}\xi (z)dE_{T}(z)$. Using this, in combination with (\ref%
{ZakEq9}), we therefore get
\begin{eqnarray*}
p(Mf_{1},Mf_{2})(z) &=&\sum_{k\in \mathbb{Z}}z^{k}\left\langle
T_{k}U^{-1}m_{0}(T)f_{1}\mid U^{-1}m_{0}(T)f_{2}\right\rangle \\
&=&\sum_{k\in \mathbb{Z}}z^{k}\left\langle T_{Nk}m_{0}(T)f_{1}\mid
m_{0}(T)f_{2}\right\rangle \\
&=&\frac{1}{N}\sum_{%
\begin{array}{c}
{\tiny w\in }\mathbb{T} \\
{\tiny w}^{{\tiny N}}{\tiny =z}%
\end{array}%
}p\left( f_{1},\left\vert m_{0}\right\vert ^{2}(T)f_{2}\right) (w)
\end{eqnarray*}%
and therefore%
\begin{eqnarray*}
&&\dint\limits_{\mathbb{T}}\xi (z)p(Mf_{1},Mf_{2})(z)d\mu (z) \\
&=&\dint\limits_{\mathbb{T}}\xi (z^{N})p(f_{1},\left\vert m_{0}\right\vert
^{2}(T)f_{2})(z)d\mu (z) \\
&=&\dint\limits_{\mathbb{T}}\xi (z^{N})\left\vert m_{0}(z)\right\vert
^{2}p(f_{1},f_{2})(z)d\mu (z) \\
&=&\dint\limits_{\mathbb{T}}\xi (z)R(p(f_{1},f_{2}))(z)d\mu (z)\text{.}
\end{eqnarray*}%
Since this is valid for all $\xi \in C(\mathbb{T})$, a comparison of the two
sides in the last formula, now yields the desired identity (\ref{ZakEq12}).
\end{proof}

\begin{remark}
\label{ZakRemarks}An immediate consequence of the lemma is that if some $%
\varphi \in \mathcal{H}$ satisfies $M\varphi =\varphi $, or equivalently%
\begin{equation}
U\varphi =\sum_{k\in \mathbb{Z}}a_{k}T_{k}\varphi  \label{ZakEq13}
\end{equation}%
then the corresponding function $h:=p(\varphi ,\varphi )$ on $\mathbb{T}$
satisfies $R(h)=h$. Recall that \emph{(}\ref{ZakEq13}\emph{)} is the \emph{%
scaling equation}. If further the functions $\{T_{k}\varphi :k\in \mathbb{Z}%
\}$ on the right hand side in (\ref{ZakEq13}) can be chosen orthogonal, then%
\begin{equation}
p(\varphi ,\varphi )(z)=\left\langle \varphi \mid \varphi \right\rangle _{%
\mathcal{H}}=\left\Vert \varphi \right\Vert _{\mathcal{H}}^{2}
\label{ZakEq14}
\end{equation}%
is the constant function, and so we are back to the special normalization
condition \emph{(}\ref{ZakEq2}\emph{)} above.
\end{remark}

The function $p(\varphi ,\varphi )$ is called the \emph{auto-correlation
function} since its Fourier coefficients%
\begin{equation}
\dint\limits_{\mathbb{T}}z^{-k}p(\varphi ,\varphi )(z)d\mu (z)=\left\langle
T_{k}\varphi \mid \varphi \right\rangle _{\mathcal{H}}  \label{ZakEq15}
\end{equation}%
are the auto-correlation numbers.

\begin{lemma}
\label{ZakLemma2}Let $m_{0}\in C(\mathbb{T})$ be given, and
suppose \textup{(}\ref{ZakEq2}\textup{)} holds.
\textup{(}a\textup{)} Then there is a probability measure $\nu
=\nu _{m_{0}}$ depending on $m_{0}$ such
that%
\begin{equation}
\dint\limits_{\mathbb{T}}\xi d\nu =\dint\limits_{\mathbb{T}}R(\xi )d\nu
\label{ZakEq16}
\end{equation}%
holds for all $\xi \in C(\mathbb{T})$.

(b) If $m_{0}$ is further assumed to be in the Lipschitz space $\limfunc{Lip}%
(\mathbb{T})$, then the following limit exists
\begin{equation}
\underset{n\rightarrow \infty }{\lim
}\frac{1}{n}\sum_{k=0}^{n-1}\mu R^{k}=\nu \label{ZakEq17}
\end{equation}%
where $\mu R^{n}(\xi ):=\mu (R^{n}\xi
)=\dint\limits_{\mathbb{T}}R^{n}\xi d\mu $, $\mu $ is the Haar
measure, and the convergence in (\ref{ZakEq17}) is in the
Hausdorff metric (details below). The measure $\nu$ satisfies
(\ref{ZakEq16}). Moreover, when then operator $R$ from
$C(\mathbb{T})$ to $C(\mathbb{T})$ has Perron-Frobenius spectrum
(i.e., $1$ is the only eigenvalue of absolute value $1$) then the
limit
\begin{equation}\label{eq3_18_1}
\lim_{n\rightarrow\infty}\mu R^n=\nu,
\end{equation}
exists and gives the unique invariant measure $\nu$.
\end{lemma}

\begin{definition}
\label{ZakDef1}A function $\xi $ on $\mathbb{T}$ is said to be in $\limfunc{%
Lip}(\mathbb{T})$ if%
\begin{equation}
D\xi :=\sup_{-\pi \leq s<t<\pi }\left\vert \frac{\xi (e^{is})-\xi (e^{it})}{%
s-t}\right\vert <\infty \text{.}  \label{ZakEq18}
\end{equation}%
The Lipschitz-norm is $\left\Vert \xi \right\Vert _{\limfunc{Lip}}:=D\xi
+\xi (1)$.
\end{definition}

The Hausdorff distance between two real valued measures $\nu _{1},\nu _{2}$
on $\mathbb{T}$ is defined as%
\begin{equation}
\mbox{dist}_{\mbox{Haus}}(\nu _{1},\nu _{2})=\sup \left\{
\dint\limits_{\mathbb{T}}\xi d\nu _{1}-\dint \xi d\nu _{2}\mid \xi
\in \limfunc{Lip}(T)\text{, }\xi \text{ real valued, and }D\xi
\leq 1\right\} \text{.}  \label{ZakEq19}
\end{equation}%
Hence the conclusion is the lemma states that if $m_{0}\in \limfunc{Lip}$
satisfies (\ref{ZakEq2}), then%
\begin{equation}
\lim_{n\rightarrow \infty }\mbox{dist}_{\mbox{Haus}}(\nu ,\mu R^{n})=0%
\text{.}  \label{ZakEq20}
\end{equation}%
Note that both the Ruelle operator $R$ and the measure $\nu $ depend on the
function $m_{0}$. Condition (\ref{ZakEq2}) states that the constant function
$\hat{1}$ is a right-Perron Frobenius eigenvector, and (\ref{ZakEq16}) that $%
\nu $ is a left-Perron-Frobenius eigenvector.\linebreak

\begin{proof}
{\it (Lemma \ref{ZakLemma2})} The lemma is essentially a special
case of the Perron-Frobenius-Ruelle theorem, see \cite{Bal00}.
Also note that an
immediate consequence of (\ref{ZakEq16}) is the invariance%
\begin{equation}
\dint\limits_{\mathbb{T}}\xi (z^{N})d\nu (z)=\dint\limits_{\mathbb{T}}\xi
(z)d\nu (z)  \label{ZakEq21}
\end{equation}%
A key step in the proof is the following estimate: Let $\left\Vert \xi
\right\Vert :=\sup_{z\in \mathbb{T}}\left\vert \xi (z)\right\vert $;\qquad
then the following estimate%
\begin{equation}
\left\Vert R^{n}\xi \right\Vert _{\limfunc{Lip}}\leq \frac{1}{N^{n}}%
\left\Vert \xi \right\Vert _{\limfunc{Lip}}+2D\left( \left\vert
m_{0}\right\vert ^{2}\right) \cdot \left\Vert \xi \right\Vert
\label{ZakEq22}
\end{equation}%
hold for all $\xi \in \limfunc{Lip}(\mathbb{T})$ and all $n\in \mathbb{Z}%
_{+} $. We leave the details of the verification of (\ref{ZakEq22}) to the
reader; see also \cite{BrJo02}, and \cite{Kea72}.
\end{proof}

For functions $m$ and $m^{\prime }$ on $\mathbb{T}$, define the form $%
\left\langle m,m^{\prime }\right\rangle $ as a function on $\mathbb{T}$ as
follows%
\begin{equation}
\left\langle m,m^{\prime }\right\rangle (z):=\frac{1}{N}\dsum\limits_{%
\underset{w^{N}=z}{{\small w\in }\mathbb{T}}}\overline{m(w)}m^{\prime }(w)%
\text{.}  \label{ZakEq23}
\end{equation}

\begin{lemma}
\label{ZakLemma3}Let $N\geq2$, and let $(\mathcal{H},U,T)$ be a
system which satisfies the commutation relation \textup{(}\ref{ZakEq0}%
\textup{)}\textit{.} Let $m$ and $m^{\prime }$ be functions on $\mathbb{T}$
which both satisfy the normalization condition \textup{(}\ref{ZakEq2}\textup{%
)}. Let $M=M_{m}$ be the cascade approximation operator, $M=U^{-1}m(T)$, and
suppose the scaling identity $M\varphi =\varphi $ has a non-zero solution in
$\mathcal{H}$ such that the vectors $T_{k}\varphi $, $k\in \mathbb{Z}$, are
mutually orthogonal. Let $M^{\prime }=M_{m^{\prime }}=U^{-1}m^{\prime }(T)$.

Then%
\begin{equation}
p(\varphi ,M^{\prime }\varphi )(z)=\left\Vert \varphi \right\Vert _{\mathcal{%
H}}^{2}\cdot \left\langle m,m^{\prime }\right\rangle (z),\qquad z\in \mathbb{%
T}.  \label{ZakEq24}
\end{equation}
\end{lemma}

\begin{proof}
Using $M\varphi =\varphi $, and $p(\varphi ,\varphi )(z)=\left\Vert \varphi
\right\Vert _{\mathcal{H}}^{2}$, we get
\begin{eqnarray*}
p(\varphi ,M^{\prime }\varphi )(z) &=&p(M\varphi ,M^{\prime }\varphi ) \\
&=&\frac{1}{N}\dsum\limits_{w^{N}=z}\overline{m(w)}m^{\prime }(w)p(\varphi
,\varphi )(w) \\
&=&\left\Vert \varphi \right\Vert _{\mathcal{H}}^{2}\cdot \frac{1}{N}%
\dsum\limits_{w^{N}=z}\overline{m(w)}m^{\prime }(w)=\left\Vert \varphi
\right\Vert _{\mathcal{H}}^{2}\cdot \left\langle m,m^{\prime }\right\rangle
(z)
\end{eqnarray*}%
which is the desired identity (\ref{ZakEq24}) in the conclusion of the lemma.
\end{proof}

\begin{corollary}
\label{ZakCorollary1}With the assumptions in Lemma \ref{ZakLemma3}, set%
\begin{equation}
m^{(n)}(z):=m(z)m(z^{N})\ldots m(z^{N^{n-1}})  \label{ZakEq25}
\end{equation}%
and%
\begin{equation*}
m^{\prime (n)}(z):=m^{\prime }(z)m^{\prime }(z^{N})\ldots m^{\prime
}(z^{N^{n-1}})\text{.}
\end{equation*}%
Then%
\begin{equation}
p(\varphi ,M^{\prime n}\varphi )(z)=\left\Vert \varphi \right\Vert _{%
\mathcal{H}}^{2}\cdot \frac{1}{N^{n}}\dsum\limits_{w^{N^{n}}=z}\overline{%
m^{(n)}(w)}m^{\prime (n)}(w)  \label{ZakEq26}
\end{equation}
\end{corollary}

\begin{proof}
A direct iteration of the argument of Lemma \ref{ZakLemma3} immediately
yields the desired identity (\ref{ZakEq26}).
\end{proof}

\begin{example}
\label{ZakEx1}Let $N=3$, and let the functions $m$ and $m^{\prime }$ be
given by $m(z)=\frac{1+z^{2}}{\sqrt{2}}$ and $m^{\prime }(z)=z^{3}m(z)$. Let%
\begin{equation}
\mathcal{H}=L^{2}(\mathcal{R},(dx)^{s})\text{, }s=\log _{3}(2)\text{,}
\label{ZakEq27}
\end{equation}%
and $\varphi =\chi _{\mathbf{C}}$; i.e., $\varphi \in \mathcal{H}$ is the
indicator function of the middle-third Cantor set, $\mathbf{C}\subset I$.
Let $Tf(x)=f(x-1)$, and $(Uf)(x)=\frac{1}{\sqrt{2}}f(\frac{x}{3})$. Then the
conditions in Lemma \ref{ZakLemma3} and Corollary \ref{ZakCorollary1} are
satisfied for this system. Specifically $\left\Vert \varphi \right\Vert _{%
\mathcal{H}}=1$, $M\varphi =\varphi $ where $M=U^{-1}m(T)$; i.e.,%
\begin{equation}
\varphi (x)=\varphi (3x)+\varphi (3x-2)\text{,}  \label{ZakEq28}
\end{equation}%
or%
\begin{equation}
3\mathbf{C}=\mathbf{C}\cup (\mathbf{C}+2)\text{.}  \label{ZakEq29}
\end{equation}%
Also notice that%
\begin{equation}
M^{\prime }=U^{-1}m^{\prime }(T)=TM\text{.}  \label{ZakEq30}
\end{equation}%
It follows that
\begin{equation*}
p(\varphi ,M^{\prime n}\varphi )(z)=z^{1+3+3^{2}+\cdots +3^{n-1}}=z^{\frac{%
3^{n}-1}{2}}\text{;}
\end{equation*}%
and therefore%
\begin{equation*}
\left\langle \varphi ,M^{\prime n}\varphi \right\rangle =\dint\limits_{%
\mathbb{T}}p(\varphi ,M^{\prime n}\varphi )(z)d\mu (z)=0
\end{equation*}%
for all $n\in \mathbb{Z}_{+}$. In contrast to the standard cascade
approximation for $L^{2}(\mathbb{R})$ with Lebesgue measure, we see that the
cascade iteration on the Cantor function $\chi _{\mathbf{C}}$; i.e.,%
\begin{equation}
\chi _{\mathbf{C}},M^{\prime }\chi _{\mathbf{C}},M^{\prime 2}\chi _{\mathbf{C%
}},\ldots  \label{ZakEq31}
\end{equation}%
does not converge in the Hilbert space $\mathcal{H}_{s}=L^{2}(\mathcal{R}%
,(dx)^{s})$. In fact the vectors in the sequence \textup{(}\ref{ZakEq31}%
\textup{)} are mutually orthogonal.
\end{example}

While the measure $\nu \in M(T)$ from Lemma \ref{ZakLemma3} is generally not
absolutely continuous with respect to the Haar measure $\mu $ on $\mathbb{T}$%
, the next result shows that it is the limit of the measures $\left\vert
m^{(n)}(z)\right\vert ^{2}d\mu (z)$ as $n\rightarrow \infty $ if the
function $m$ is given to satisfy (\ref{ZakEq2}) and if the sequence $m^{(n)}$
is defined by (\ref{ZakEq25}).

\begin{proposition}
\label{ZakProp1}Let $m\in L^{\infty }(\mathbb{T})$ be given. Suppose \textup{%
(}\ref{ZakEq2}\textup{)} holds, and let $\nu $ be the Perron-Frobenius
measure of Lemma \ref{ZakLemma3}; i.e., the measure $\nu =\nu _{m}$ arising
as a limit \textup{(}\ref{eq3_18_1}\textup{)}. Then%
\begin{equation}
\lim_{n\rightarrow \infty }\dint\limits_{\mathbb{T}}\xi (z)\left\vert
m^{(n)}(z)\right\vert ^{2}d\mu (z)=\dint\limits_{\mathbb{T}}\xi (z)d\nu (z)
\label{ZakEq32}
\end{equation}%
holds for all $\xi \in C(\mathbb{T})$.
\end{proposition}

\begin{proof}
Calculating the integrals on the left-hand side in (\ref{ZakEq32}), we get
\begin{equation*}
\dint\limits_{\mathbb{T}}\xi \left\vert m^{(n)}\right\vert ^{2}d\mu
=\dint\limits_{\mathbb{T}}\xi R^{\ast ^{n}}(\hat{1})d\mu =\dint\limits_{%
\mathbb{T}}R^{n}(\xi )d\mu \underset{n\rightarrow \infty }{\longrightarrow }%
\dint\limits_{\mathbb{T}}\xi d\nu
\end{equation*}%
where (\ref{ZakEq17}) was used in the last step. This is the desired
conclusion (\ref{ZakEq32}) of the proposition.
\end{proof}

\begin{corollary}
\label{ZakCorollary2}If $m(z)=\frac{1+z^{2}}{\sqrt{2}}$ is the function from
Example \ref{ZakEx1}, then the substitution $z=e^{it}$ yields the following
limit formula for the corresponding Perron-Frobenius measure $\nu $, written
in multiplicative notation:%
\begin{equation}
\lim_{n\rightarrow \infty }\frac{1}{2\pi }\dprod\limits_{k=1}^{n}\left(
1+\cos (2\cdot 3^{k}t)\right) =d\nu (t)  \label{ZakEq33}
\end{equation}
\end{corollary}

\begin{proof}
The expression on the left-hand side in (\ref{ZakEq33}) is called
a Riesz-product, and it belongs to a wider family of examples; see
for example \cite{Kea72}, and \cite{Zyg68}. The limit measure $\nu
$ is known to be singular. It follows, for example from
\cite{Zyg68}. A computation shows
that the Fourier coefficients $\hat{\nu}(n):=\int_{\mathbb{T}}z^{n}d\nu (z)$%
, are real valued, satisfy $\hat{\nu}(0)=1$, $\hat{\nu}(1)=0$, $\hat{\nu}%
(-n)=\hat{\nu}(n)$, $\hat{\nu}(3n)=\hat{\nu}(n)$, $\hat{\nu}(3n-2)=\frac{1}{2%
}\hat{\nu}(n)$, $\hat{\nu}(3n+2)=\frac{1}{2}\hat{\nu}(n)$, for all
$n\in
\mathbb{Z}$. Hence
\begin{equation}\label{eqwiener}
\lim \frac{1}{2n+1}\sum_{k=-n}^{n}\left\vert \hat{\nu}(k)\right\vert ^{2}=0=\mbox{sum of atoms};
\end{equation}
i .e., $\sum \left\vert \nu(\{z\})\right\vert ^{2}=0$.
The last conclusion is from Wiener's theorem, and implies that $\nu $ has no
atoms; i.e., $\nu (\{z\})=0$ for all $z\in \mathbb{T}$.
\par
To check (\ref{eqwiener}) some computations are required.
Let
$$s_n=|\hat\nu(0)|^2+|\hat\nu(1)|^2+...+|\hat\nu(n)|^2,\quad(n\in\mathbb{N}).$$
Then, using the recursive relations for $\hat\nu$ we have
\begin{align*}
s_{3^{n+1}}&=s_{3^n}+\sum_{k=3^n+1,k\equiv 0\mod 3}^{3^{n+1}}|\hat\nu(k)|^2\\
&+\sum_{k=3^n+1,k\equiv 2\mod 3}^{3^{n+1}}|\hat\nu(k)|^2+\sum_{k=3^n+1,k\equiv -2\mod 3}^{3^{n+1}}|\hat\nu(k)|^2\\
&\leq s_{3^n}+s_{3^n}+\frac14 s_{3^n}+\frac14 s_{3^n}=\frac52 s_{3^n}.
\end{align*}
By induction
$$s_{3^n}\leq \left(\frac52\right)^n s_0.$$
Now take $k$ arbitrary then for some $n$, $k$ is inbetween $3^n$ and $3^{n+1}$ so
$$\frac{s_k}{k}\leq\frac{s_{3^{n+1}}}{3^n}\leq\frac{\left(\frac52\right)^{n+1}s_0}{3^n}=\left(\frac56\right)^n\frac52 s_0,$$
which shows that $s_k/k$ converges to $0$ and this proves (\ref{eqwiener}).
\par
In summary, the
measure $\nu $ is singular and non-atomic. In the next section we show that $%
\nu $ has full support.

Moreover there are supporting sets for $\nu $ which have zero Haar measure
as subsets of $\mathbb{T}$. Concrete constructions are given below.

It follows from the recursive relations for the numbers $\hat{\nu}(n)$ that $%
\hat{\nu}(2k+1)=0$, $k\in \mathbb{Z}$; i.e., that all the odd Fourier
coefficients vanish.

Each integer $n\in\mathbb{Z}_+$ has a representation of the
following form:
\begin{equation}
n=l_{0}+l_{1}\cdot 3+l_{2}\cdot 3^{2}+\cdots +l_{p}\cdot 3^{p}\text{,\qquad }%
l_{i}\in \{0,-2,2\}, i<p, l_p\in\{1,2\}.  \label{ZakEq34}
\end{equation}%
The representation is not unique (because for example $1=3\cdot
1-2$) but uniqueness is obtain if we impose that $l_{p-1}\neq -2$
when $l_p=1$. We define a counting function $\#(n)$ which records
the occurrence of values $-2$ and $2$ for the `trigets' $l_{i}$.
Hence, if $n\in
\mathbb{Z}_{+}$ is even, then%
\begin{equation*}
\hat{\nu}(n)=2^{-\#(n)}\text{.}
\end{equation*}

Now, introduce the following sequence of
functions%
\begin{equation}
g_{k}(z):=z^{2\cdot 3^{k}}-\frac{1}{2}\text{,\qquad }z\in \mathbb{T}\text{.}
\label{ZakEq35}
\end{equation}%
with inner products as follows with respect to the measure $\nu $:%
\begin{equation}
\left\langle g_{k}\mid g_{l}\right\rangle _{\nu }=\dint\limits_{\mathbb{T}}%
\overline{g_{k}(z)}g_{l}(z)d\nu (z)=\frac{3}{4}\delta _{k,l}  \label{ZakEq36}
\end{equation}%
It follows from this that the series%
\begin{equation}
g(z):=\dsum\limits_{k=0}^{\infty }\frac{1}{k+1}g_{k}(z)\text{,\qquad }z\in
\mathbb{T}\text{,}  \label{ZakEq37}
\end{equation}%
is convergent in $L^{2}(\mathbb{T},\nu )$ with%
\begin{equation}
\left\Vert g\right\Vert _{\nu }^{2}=\dsum\limits_{k=0}^{\infty }\frac{1}{%
(k+1)^{2}}\left\Vert g_{k}\right\Vert _{\nu }^{2}=\frac{\pi ^{2}}{8}\text{.}
\label{ZakEq38}
\end{equation}

Using the Riesz-Fisher theorem, we get a Borel subset $A\subset \mathbb{T}$,
$\nu (A)=1$; i.e., $\nu (\mathbb{T\diagdown }A)=0$, and a subsequence $%
n_{1}<n_{2}<n_{3}<\ldots ,n_{i}\rightarrow \infty $, such that the series%
\begin{equation}
\dsum\limits_{k=0}^{n_{i}}\frac{1}{k+1}g_{k}(z)  \label{ZakEq39}
\end{equation}%
is pointwise convergent, $i\rightarrow \infty $, for all $z\in A$. But note
that%
\begin{equation}
\dsum\limits_{k}\frac{1}{k+1}g_{k}(z)=\dsum\limits_{k}\frac{1}{k+1}z^{2\cdot
3^{k}}-\frac{1}{2}\dsum\limits_{k}\frac{1}{k+1}\text{.}  \label{ZakEq40}
\end{equation}%
Using now Carleson's theorem about Fourier series on $\mathbb{T}$ with
respect to Haar measure $\mu $ (=Lebesgue measure), we conclude that there
is a Borel subset $B\subset \mathbb{T}$, $\mu (B)=1$; i.e., $\mu (\mathbb{%
T\diagdown }B)=0$, such that the series $\sum_{k=0}^{n}\frac{1}{k+1}%
z^{2\cdot 3^{k}}$ is pointwise convergent, $n\rightarrow \infty $, for all $%
z\in B$. But identity (\ref{ZakEq40}) implies that $A\cap B=\varnothing $,
and so $\mu (A)=0$. The supporting set $A$ for the measure $\nu $ has
Lebesgue measure zero; and moreover the two measures $\nu $ and $\mu $
(=Haar measure on $\mathbb{T}$) are mutually singular.
\end{proof}

It can be shown, using a theorem of Nussbaum \cite{Nus98}, that the Ruelle
operator
\begin{equation*}
(Rf)(z)=\frac{1}{3}\dsum\limits_{w^{3}=z}\left\vert m(w)\right\vert ^{2}f(w)=%
\frac{1}{3}\dsum\limits_{w^{3}=z}\left( 1+\frac{w^{2}+w^{-2}}{2}\right) f(w)
\end{equation*}%
has Perron-Frobenius spectrum on $C(\mathbb{T)}$, specifically that
\begin{equation*}
E_{R}(1)=\left\{ f\in C(\mathbb{T})\mid Rf=f\right\} =\mathbb{C}\hat{1}\text{%
,}
\end{equation*}%
and if $\lambda \in \mathbb{C}$, $\left\vert \lambda \right\vert =1$, and $%
\lambda \neq 1$, that
\begin{equation*}
E_{R}(\lambda )=\left\{ f\in C(\mathbb{T})\mid Rf=\lambda f\right\} =0\text{.%
}
\end{equation*}%
As a consequence we get that
\begin{equation*}
\left\Vert R^{n}(f)-\nu (f)\right\Vert _{\limfunc{Lip}}\leq \frac{1}{3^{n}}%
\left\Vert f\right\Vert _{\limfunc{Lip}}
\end{equation*}%
holds for all $f\in \limfunc{Lip}(\mathbb{T})$ where $\left\Vert \cdot
\right\Vert _{\limfunc{Lip}}$ denotes the Lipschitz norm on functions on $%
\mathbb{T}$.

\section{\label{SupportP&F}The support of Perron-Frobenius measures}

In this section we consider the support of the measure $\nu $. Caution: By
the support of $\nu $, we mean the support of $\nu $ when it is viewed as a
distribution; i.e., the support of $\nu $ is the complement of the union of
the open subsets in $\mathbb{T}$ where $\nu $ acts as the zero distribution,
or the zero Radon measure. Even though the support of $\nu $ may be all of $%
\mathbb{T}$, there can still be Borel subsets $E\subset \mathbb{T}$ with $%
\nu (E)=1$, but $E$ having zero Lebesgue measure.

\begin{theorem}
\label{PerFrobTheo1}Let $m_{0}\in \limfunc{Lip}(\mathbb{T})$. Suppose
\textup{(}\ref{ZakEq2}\textup{)} holds and $m_{0}$ has finitely many zeros
and let $\nu $ be the Perron-Frobenius measure of Lemma \ref{ZakLemma2}. Then exactly one of the following affirmations is
true:

\noindent \textup{(}i\textup{)} The support of $\nu $ is $\mathbb{T}$.

\noindent \textup{(}ii\textup{)} $\nu $ is atomic and the support of $\nu $
is a union of cycles $C=\{z_{1},\ldots ,z_{p}\}$ with $z_{1}^{N}=z_{2},%
\ldots ,z_{p-1}^{N}=z_{p}$, $z_{p}^{N}=z_{1}$, and $\left\vert
m_{0}\right\vert ^{2}(z_{i})=N$ for all $i\in \{1,\ldots ,p\}$.
(Such cycles are called $(m_0,N)$-cycles).
\end{theorem}

\begin{proof}
To simplify the notation, we use%
\begin{equation*}
W:=\left\vert m_{0}\right\vert ^{2}\text{.}
\end{equation*}

Note that proposition \ref{ZakProp1} gives us $\nu $ as an infinite
product $\dprod\limits_{k=1}^{\infty }W(z^{N^{k}})d\mu $. In the next lemma
we analyze the measures given by the tails of this product.

\begin{lemma}
\label{PerFrobLem1}Fix $n\geq 0$. We then have$:$

\noindent \textup{(i)} For all $\QTR{up}{f}\in C(\mathbb{T})$ the following
limit exists, and defines a measure on $\mathbb{T}\emph{:}$ $\nu
_{n}(f):=\lim_{k\rightarrow \infty }\int_{\mathbb{T}}W(z^{N^{n}})\cdots
W(z^{N^{n+k}})f(z)d\mu $

\noindent \textup{(ii)}$\nu _{n}(f)=\nu (R_{1}^{n}f)$,$ (f\in C(%
\mathbb{T)}$ where $R_{1}$ is given as in \textup{(\ref{ZakProp1})}
but with $m_{0}=1$.

\noindent \textup{(iii) }$\int_{\mathbb{T}}f(z)W(z)\cdots W(z^{N^{n-1}})d\nu
_{n}=\int_{\mathbb{T}}f(z)d\nu $, $f\in C(\mathbb{T)}$. \textup{(So }$%
\nu $ is absolutely continuous with respect to $\nu _{n}$.\textup{)}

\noindent \textup{(iv) }$\lim_{n\rightarrow \infty }\nu _{n}(f)=\mu (f)$%
, $f\in C(\mathbb{T)}$.
\end{lemma}

\begin{proof}[Proof of Lemma]
We can use a change of variable to compute
\begin{eqnarray*}
&&\lim_{k\rightarrow \infty }\int_{\mathbb{T}}W\left( z^{N^{n}}\right)
\cdots W\left( z^{N^{n+k}}\right) f(z)d\mu \\
&=&\lim_{k\rightarrow \infty }\int_{\mathbb{T}}W(z)\cdots W\left(
z^{N^{k}}\right) R_1^{n}f(z)d\mu \\
&=&\nu \left( R_{1}^{n}f\right) .
\end{eqnarray*}

This proves (i) and (ii). (iii) is immediate from proposition \ref{ZakProp1}.

For (iv) note that for $f\in C(\mathbb{T)}$ and $n\in \mathbb{N}$.%
\begin{equation*}
R_{1}^{n}f(\theta )=\frac{1}{N^{n}}\dsum\limits_{k=0}^{N^{n}-1}f\left( \frac{%
\theta +2k\pi }{N^{n}}\right)
\end{equation*}%
therefore $R_{1}^{n}f$ converges uniformly to $\int_{\mathbb{T}}fd\mu $.

Then, with (ii), $\lim_{n\rightarrow \infty }\nu _{n}(f)=\mu (f)$.
\end{proof}

We continue now the proof of the theorem. We distinguish two cases:

\noindent \textbf{Case I:} All measures $\nu _{n}$ are absolutely continuous
with respect to $\nu $. In this case we prove that the support of $\nu $ is $%
\mathbb{T}$. Assume the contrary. Then there is an open set $U$ with $\nu
(U)=0$. This implies $\nu _{n}(U)=0$ for all $n$.

Take $f\in C(\mathbb{T)}$ with support contained in $U$. Then $\nu
_{n}(f)=0 $. Take the limit and use lemma \ref{PerFrobLem1}(iv),
it follows
that $\mu (f)=0$. As $f$ is arbitrary, $\mu (U)=0$. But this implies $%
U=\varnothing $, so the support of $\nu $ is indeed $\mathbb{T}$.

\noindent \textbf{Case II:} There is an $n\in \mathbb{N}$ such that $\nu
_{n} $ is not absolutely continuous with respect to $\nu $. This means that
there is a Borel set $E$ with $\nu (E)=0$ and $\nu _{n}(E)>0$.

We prove that there is a $\limfunc{zero}$ of $W^{(n)}$, call it
$z_{0}$, such that $\nu _{n}(\{z_{0}\})>0$. Suppose not. Then,
take $E^{\prime }=E\diagdown \limfunc{zeros}(W^{(n)})$, $\nu
(E^{\prime })=0$ $\nu _{n}(E^{\prime })>0.$

The measure $\nu $ is regular so there is a compact subset $K$ of $E^{\prime
}$ such that $\nu _{n}(K)>0$. Of course $\nu (K)=0$. Since $K$ has no $%
\limfunc{zeros}$ of $W^{(n)}$ and this is continuous, $W^{(n)}$ is bounded
away from $0$ on $K$. Then, with lemma \ref{PerFrobLem1}(iii)%
\begin{equation*}
0=\dint\limits_{K}\frac{1}{W^{(n)}(z)}d\nu =\dint\limits_{K}W^{(n)}(z)\frac{1}{%
W^{(n)}(z)}d\nu _{n}=\nu _{n}(K)
\end{equation*}%
which is a contradiction.

Thus, there is a $z_{0}\in \limfunc{zeros}(W^{(n)})$ with $\nu (\{z_{0}\})>0$%
.

We know also that $\nu (f)=\nu (Rf)$ for all $f\in C(\mathbb{T)}$. By
approximation (Lusin's theorem) the same equality is true for all bounded
Borel functions. Then
\begin{equation}\label{eq3_41_1}
0<\nu (\chi _{\{z_{0}\}})=\nu (R\chi _{\{z_{0}\}})=\nu \left( \frac{1}{N}%
\dsum\limits_{w^{N}=z}W(w)\chi _{\{z_{0}\}}(w)\right)
=\frac{W(z_{0})}{N}\nu (\chi _{\{z_{0}^{N}\}})\text{.}
\end{equation}
Therefore $W(z_{0})>0$ and $\nu (\{z_{0}^{N}\})>0$. By induction, $%
W(z_{0}^{N^{k}})>0$ and $\nu (\{z_{0}^{N^{k}}\})$ for all $k\in \mathbb{N}$.

Since (\ref{ZakEq2}) holds, $W(z)\leq N$ for all $z\in
\mathbb{T}$; so from the previous computation we obtain
\begin{equation*}
\nu \left( \left\{ z^{N}\right\} \right) \geq \nu (\{z\})\text{.}
\end{equation*}

Since $\nu $ is a finite measure, the orbit $\left\{ z_{0}^{N^{k}}\mid k\in
\mathbb{N}\right\} $ has to be finite so the points $z_{0},z_{0}^{N},\ldots
,z_{0}^{N^{p}}$ will form a cycle, $z_{0}^{N^{p+1}}=z_{0}$. Also,
\begin{equation*}
\nu (\{z_{0}\})=\nu \left( \left\{ z_{0}^{N^{p+1}}\right\} \right) \geq \nu
\left( \left\{ z_{0}^{N^{p}}\right\} \right) \geq \cdots \geq \nu \left(
\left\{ z_{0}\right\} \right)
\end{equation*}%
hence all inequalitites are in fact equalities and with (\ref{eq3_41_1}), this shows that $%
W(z_{0}^{N^{k}})=N$ for $k\in \{0,\ldots ,p\}$.

We are now in the \textquotedblleft classical\textquotedblright\ case and we
can use corollary 2.18 in \cite{Dut1} (See also \cite{BrJo02}) to conclude
that $\nu $ must be atomic and supported on cycles as mentioned in the
theorem.
\end{proof}
\begin{lemma}
\label{PerFrobLem2}Let $m_{0},m_{0}^{\prime }$ be Lipschitz, with finitely
many zeros. Suppose
\begin{equation*}
R_{m_{0}}(\hat{1})=\hat{1}=R_{m_{0}^{\prime }}(\hat{1})\text{,}
\end{equation*}
and suppose there are no $m_{0}$ or $m_{0}^{\prime }$-cycles. Assume in
addition that $m_{0}$ and $m_{0}^{\prime }$ have the same Perron-Frobenius
measure $\nu $, then $\left\vert m_{0}\right\vert =\left\vert m_{0}^{\prime
}\right\vert $.
\end{lemma}

\begin{proof}
With $W:=\left\vert m_{0}\right\vert ^{2}$, we have, from proposition \ref{ZakProp1}, for $f\in C(\mathbb{T)}$%
:
\begin{equation*}
\lim_{n\rightarrow \infty }\dint\limits_{\mathbb{T}}f(z)W^{(n)}(z)d\mu
(z)=\dint\limits_{\mathbb{T}}f(z)d\nu (z)\text{.}
\end{equation*}%
Then
\begin{eqnarray*}
\lim_{n\rightarrow \infty }\dint\limits_{\mathbb{T}}f(z)W^{(n)}(z^{N})d\mu
&=&\lim_{n\rightarrow \infty }\dint\limits_{\mathbb{T}}R_{1}f(z)W^{(n)}(z)d%
\mu \\
&=&\dint\limits_{\mathbb{T}}R_{1}f(z)d\nu \text{.}
\end{eqnarray*}%
Take $f\in C(\mathbb{T)}$ such that $f$ is zero in a neighborhood of $%
\limfunc{zeros}(m_{0})$. Then $\frac{f}{\left\vert m_{0}\right\vert ^{2}}$
is continuous. So
\begin{equation*}
\lim_{n\rightarrow \infty }\dint\limits_{\mathbb{T}}\frac{f(z)}{\left\vert
m_{0}(z)\right\vert ^{2}}W^{(n)}(z)d\mu (z)=\dint\limits_{\mathbb{T}}\frac{%
f(z)}{\left\vert m_{0}(z)\right\vert ^{2}}d\nu (z)\text{.}
\end{equation*}%
On the other hand%
\begin{eqnarray*}
\lim_{n\rightarrow \infty }\dint\limits_{\mathbb{T}}\frac{f(z)}{\left\vert
m_{0}(z)\right\vert ^{2}}W^{(n)}(z)d\mu (z) &=&\lim_{n\rightarrow \infty
}\dint\limits_{\mathbb{T}}\frac{f(z)}{\left\vert m_{0}(z)\right\vert ^{2}}%
\left\vert m_{0}(z)\right\vert ^{2}\cdots \left\vert
m_{0}(z^{N^{n-1}})\right\vert ^{2}d\mu \\
&=&\lim_{n\rightarrow \infty }\dint\limits_{\mathbb{T}}f(z)W^{(n-1)}(z^{N})d%
\nu (z) \\
&=&\dint\limits_{\mathbb{T}}R_{1}f(z)d\nu \text{.}
\end{eqnarray*}%
Thus
\begin{equation}
\dint\limits_{\mathbb{T}}\frac{f(z)}{\left\vert m_{0}(z)\right\vert ^{2}}%
d\nu (z)=\dint\limits_{\mathbb{T}}R_{1}f(z)d\nu (z)  \label{PerFrobEq1}
\end{equation}%
for all $f\in C(\mathbb{T)}$ which are zero in a neighborhood of $\limfunc{%
zeros}(m_{0})$.

The same argument can be applied to $m_{0}^{\prime }$. But note that the
right-hand side of (\ref{PerFrobEq1}) doesn't depend on $m_{0}$ or $%
m_{0}^{\prime }$ (because $\nu $ is the same).

Therefore
\begin{equation}
\dint\limits_{\mathbb{T}}\frac{f(z)}{\left\vert m_{0}(z)\right\vert ^{2}}%
d\nu (z)=\dint\limits_{\mathbb{T}}\frac{f(z)}{\left\vert m_{0}^{\prime
}(z)\right\vert ^{2}}d\nu (z)  \label{PerFrobEq2}
\end{equation}%
for all $f\in C(\mathbb{T)}$ which are zero on a neighborhood of $\limfunc{%
zeros}(m_{0})\cup \limfunc{zeros}(m_{0}^{\prime })$.

From (\ref{PerFrobEq2}) it follows that $\left\vert
m_{0}\right\vert =\left\vert m_{0}^{\prime }\right\vert $, $\nu
$-almost everywhere on
\begin{equation*}
C:=\mathbb{T}\diagdown (\limfunc{zeros}(m_{0})\cup \limfunc{zeros}%
(m_{0}^{\prime }))\text{.}
\end{equation*}
Since the support of $\nu $ is $\mathbb{T}$, this implies that $\left\vert
m_{0}\right\vert =\left\vert m_{0}^{\prime }\right\vert $ on a set which is
dense in $\mathbb{T}$, and since the zeros of $m_{0}$ and $m_{0}^{\prime }$ are
finite in number, $\left\vert m_{0}\right\vert =\left\vert m_{0}^{\prime
}\right\vert $ on a dense subset of $\mathbb{T}$. By continuity therefore
\begin{equation*}
\left\vert m_{0}\right\vert =\left\vert m_{0}^{\prime }\right\vert \text{ on
}\mathbb{T}\text{.}
\end{equation*}
\end{proof}

If $m_{0}$ and $W:=\left\vert m_{0}\right\vert ^{2}$ are not assumed
continuous, there is still a variant of Theorem \ref{PerFrobTheo1}, but with
a weaker conclusion. For functions $W$ on $\mathbb{T}$ we introduce the
following axioms. The a.e. conditions are taken with respect to the Haar
measure $\mu $ on $\mathbb{T}$:

\noindent \textup{(i) }$W\in L^{\infty }(\mathbb{T)}$,

\noindent \textup{(ii) }$W\geq 0$ a.e. on $\mathbb{T}$ with respect to $\mu $%
,

\noindent \textup{(iii) }$\frac{1}{N}\dsum\limits_{w^{N}=z}W(w)=1$ a.e. on $%
\mathbb{T}$,

\noindent \textup{(iv) The limit }%
\begin{equation}
d\nu _{W}:=\lim_{n\rightarrow \infty }\frac{1}{n}\dsum%
\limits_{k=1}^{n}W^{(k)}d\mu  \label{PerFrobEq3}
\end{equation}%
exists in $M_{1}(\mathbb{T)}$. Here, as before,
\begin{equation}
W^{(k)}(z):=W(z)W(z^{N})\cdots W\left( z^{N^{k-1}}\right) \text{,}
\label{PerFrobEq4}
\end{equation}%
and

\begin{equation}
(R_{W}f)(z)=\frac{1}{N}\dsum\limits_{w^{N}=z}W(w)f(w)\text{,}\qquad z\in
\mathbb{T}\text{,\qquad }f\in C(\mathbb{T)}\text{.}  \label{PerFrobEq5}
\end{equation}

\begin{theorem}
\label{PerFrobTheo2}Let $N\in \mathbb{Z}_{+}$, and let a function $W$ be
given on $\mathbb{T}$ satisfying \textup{(i)--(iv) }above. Let $\nu _{W}$, $%
W^{(k)}$ and $R_{W}$ be given as in \textup{(\ref{PerFrobEq3})--(\ref%
{PerFrobEq5}).} Then \textup{(}a\textup{)}--\textup{(}c\textup{)} hold\emph{:%
}

\noindent \textup{(}a\textup{)} $\nu _{W}\left( R_{W}(f)\right) =\nu _{W}(f)$%
,\qquad $f\in C(\mathbb{T)}$.

\noindent \textup{(}b\textup{)} If $W_{0}$, and $W$ both satisfy \textup{(}i%
\textup{)}--\textup{(}iv\textup{)}, set
\begin{equation}
d\nu _{0}:=\lim_{n\rightarrow \infty }\frac{1}{n}\dsum%
\limits_{k=1}^{n}W_{0}^{(k)}d\mu \text{,}  \label{PerFrobEq6}
\end{equation}%
and
\begin{equation*}
f\longrightarrow \nu _{0}\left( R_{W}(fW_{0})\right)
\end{equation*}%
defines a measure on $\mathbb{T}$ which is absolutely continuous with
respect to $\nu _{0}$ with Radon-Nikodym derivative $W$\emph{;} i.e., we
have
\begin{equation}
\nu _{0}\left( R_{W}(fW_{0})\right) =\nu _{0}(fW)\text{,\qquad }f\in C(%
\mathbb{T)}\text{.}  \label{PerFrobEq7}
\end{equation}

\noindent \textup{(}c\textup{)} If $W_{0}$ and $W$ are as in \textup{(}b%
\textup{)}, then the following are equivalent\emph{:}

\noindent \textup{(}c$_{\mathit{1}}$\textup{)} $\nu _{0}R_{W}=\nu _{0}$

and

\noindent \textup{(}c$_{\mathit{2}}$\textup{)} There is a Borel subset $%
E\subset \mathbb{T}$ such that $\nu _{0}(E)=1$ and
\begin{equation}
W_{0}(z)=W(z)  \label{PerFrobEq8}
\end{equation}%
for all $z\in E$.
\end{theorem}

\begin{proof}
The structure of this proof is as that of Theorem \ref{PerFrobTheo1}, the
essential step consists of the following two duality identities. Each one
amounts to a basic property of the Haar measure $\mu $ on $\mathbb{T}$. For
every $k\in \mathbb{Z}_{+}$ and $f\in C(\mathbb{T)}$, we have
\begin{equation*}
\dint\limits_{\mathbb{T}}R_{W}(f)W^{(k)}d\mu =\dint\limits_{\mathbb{T}%
}fW^{(k+1)}d\mu
\end{equation*}%
and
\begin{equation*}
\dint\limits_{\mathbb{T}}R_{W}(fW_{0})W_{0}^{(k)}d\mu =\dint\limits_{\mathbb{%
T}}fWW_{0}^{(k+1)}d\mu \text{.}
\end{equation*}%
Using these, and \textup{(\ref{PerFrobEq1}) for the measure }$\nu _{0}$, the
desired conclusions follow as before. For \textup{(}c\textup{)}, in
particular, we note that \textup{(}c$_{1}$\textup{)} yields the identity $%
Wd\nu _{0}=W_{0}d\nu _{0}$ for the two functions $W_{0}$ and $W$ on $\mathbb{%
T}$. Hence $W_{0}$ and $W$ must agree on a Borel subset in $\mathbb{T}$ of
full $\nu _{0}$-measure, and conversely.
\end{proof}

\section{\label{XformRules}Transformation Rules}

If $m_{0}$:~$\mathbb{T}\rightarrow \mathbb{C}$ is a Fourier polynomial;
i.e., represented by a finite sum $m_{0}(z)=\sum_{k}a_{k}z^{k}$, and if
\begin{equation}
\frac{1}{N}\dsum\limits_{w^{N}=z}\left\vert m_{0}(w)\right\vert ^{2}=1\text{%
,\qquad }z\in \mathbb{T}\text{,}  \label{XformEq1}
\end{equation}%
for some $N\in \mathbb{Z}_{+}$, $N\geq2$ we showed in
\cite{BrJo01a} that there are functions $m_{1},\ldots ,m_{N-1}$ on
$\mathbb{T}$ such that
\begin{equation}
\frac{1}{N}\dsum\limits_{w^{N}=z}\overline{m_{j}(w)}m_{k}(w)=\delta _{j,k}%
\text{,\qquad }z\in \mathbb{T}\text{,}  \label{XformEq2}
\end{equation}%
or equivalently, the $N\times N$ matrix
\begin{equation}
\frac{1}{\sqrt{N}}\left( m_{j}\left( e^{i\frac{k\cdot 2\pi }{N}}z\right)
\right) _{j,k=0}^{N-1}\text{,\qquad }z\in \mathbb{T}  \label{XformEq3}
\end{equation}%
is unitary; i.e., defines a function from $\mathbb{T}$ into the group $U_{N}(%
\mathbb{C})$ of all unitary $N\times N$ matrices. Moreover the functions $%
m_{1},\ldots ,m_{N-1}$ may be chosen to be Fourier polynomials of the same
total degree as $m_{0}$. For the example in Section 2, $N=3$, and $m_{0}(z)=%
\frac{1+z^{2}}{\sqrt{2}}$, condition (\ref{XformEq1}) is satisfied, and the
other two functions may be chosen as: $m_{1}(z)=z$, and $m_{2}(z)=\frac{%
1-z^{2}}{\sqrt{2}}$. There is a more general result, also from \cite{BrJo01a}
which defines a transitive action of the group $G_{N}$ of all unitary matrix
functions ($G_{N}$ is often called the $N$th order loop-group, and it is
used in homotopy theory). An element in $G_{N}$ is a function $A$:~$\mathbb{T%
}\rightarrow U_{N}(\mathbb{C})$. Let $F_{N}$ denote all the functions $%
m=(m_{j})_{j=0}^{N-1}$:~$\mathbb{T\rightarrow C}^{N}$ which satisfy (\ref%
{XformEq2}), or equivalently (\ref{XformEq3}). Define the action of $A$ on $%
m $ as follows:
\begin{equation}
m^{A}(z):=A(z^{N})m(z)\text{,\qquad }z\in \mathbb{T}  \label{XformEq4}
\end{equation}

\begin{lemma}
\label{XformLem1}The action of $G_{N}$ on $F_{N}$ is transitive and
effective; specifically, for any two $m$ and $m^{\prime }\in F_{N}$%
, there is a unique $A\in G_{N}$ such that $m^{\prime }=m^{A}$; i.e., $A$
transforms $m$ to $m^{\prime }$.
\end{lemma}

\begin{proof}
For functions $f$ and $g$ on $\mathbb{T}$,
\begin{equation}
\left\langle f,g\right\rangle _{N}(z):=\frac{1}{N}\dsum\limits_{w^{N}=z}%
\overline{f(w)}g(w)\text{,\qquad }z\in \mathbb{T}\text{.}  \label{XformEq5}
\end{equation}%
If $m,m^{\prime }\in F_{N}$ are given, set
\begin{equation}
A_{j,k}(z):=\left\langle m_{k},m_{j}\right\rangle _{N}(z)\text{,\qquad }z\in
\mathbb{T}\text{.}  \label{XformEq6}
\end{equation}%
Then an easy verification shows that $A=\left( A_{j,k}\right) _{j,k=0}^{N-1}$
defines an element in $G_{N}$ which transforms $m$ to $m^{\prime }$.
Conversely, if $m\in F_{N}$ is given, and $m^{A}$ is defined by (\ref%
{XformEq4}), then it follows that $m^{A}\in F_{N}$ if and only if $A\in
G_{N} $.
\end{proof}

If a function $m_{0}$:~$\mathbb{T}\rightarrow \mathbb{C}$ is given, and
satisfies (\ref{XformEq1}) for some $N$, then the Ruelle transfer operator
\begin{equation}
\left( Rf\right) (z):=\frac{1}{N}\dsum\limits_{w^{N}=z}\left\vert
m_{0}(w)\right\vert ^{2}f(w)\text{,\qquad }z\in \mathbb{T}  \label{XformEq7}
\end{equation}%
satisfies $R\hat{1}=\hat{1}$ where $\hat{1}$ denotes the constant function $%
1 $ on $\mathbb{T}$.

By a \textit{probability measure} on $\mathbb{T}$, we mean a (positive)
Borel measure $\nu $ on $\mathbb{T}$ such that $\nu (\mathbb{T)=}1$. The
probability measure will be denoted $M_{1}(\mathbb{T)}$.

\textbf{Terminology for measures }$\nu $ \textbf{on }$\mathbb{T}$: If $f\in
C(\mathbb{T)}$, set%
\begin{equation*}
\nu (f):=\dint\limits_{\mathbb{T}}fd\nu =\dint\limits_{\mathbb{T}}f(z)d\nu
(z)\text{.}
\end{equation*}

If $\tau $:~$\mathbb{T\rightarrow T}$ is a measurable transformation, set
\begin{equation*}
\nu ^{\tau ^{-1}}(E)=\nu (\tau ^{-1}(E))
\end{equation*}%
for Borel sets $E\subset \mathbb{T}$.

If $R$:~$C(\mathbb{T)\rightarrow }C(\mathbb{T)}$ is linear, set
\begin{equation*}
(\nu R)(f)=\nu (Rf)=\dint\limits_{\mathbb{T}}(Rf)d\nu \text{.}
\end{equation*}

If $m_{0}$ is given as in (\ref{HausEq7}), we introduce
\begin{equation}
\mathcal{L}(m_{0}):=\left\{ \nu \in M_{1}(\mathbb{T)}:\nu R_{m_{0}}=\nu
\right\} \text{.}  \label{XformEq8}
\end{equation}

If $\nu \in \mathcal{L}(m_{0})$ and $E\subset \mathbb{T}$ is a
Borel subset, we recall that $E$ \textit{supports} $\nu $ if $\nu
(E)=1$. Note that from the
examples in Section 4, it may be that the support of $\nu $ is all of $%
\mathbb{T}$ even though $\nu $ has supporting Borel sets $E$ with zero Haar
measure; i.e., $\nu (E)=1$ and $\mu (E)=0$.

We now return to the case of the middle-third Cantor set $\mathbf{C}$. Set $%
s:=\log _{3}(2)$, and view $\chi _{\mathbf{C}}$ as an element in the Hilbert
space $L^{2}(\mathbb{R},(dx)^{s})$. We recall the usual unitary operators $%
(Uf)(x):=\frac{1}{\sqrt{2}}f(\frac{x}{3})$, and $(T_{k}f)(x)=f(x-k)$, $k\in
\mathbb{Z}$, and the relation
\begin{equation}
UT_{k}U^{-1}=T_{3k}\text{,\qquad }k\in \mathbb{Z}\text{.}  \label{XformEq9}
\end{equation}

If $m_{0}(z)=\sum_{k}a_{k}z^{k}$, $m_{0}\in L^{\infty }(\mathbb{T)}$, is
given, we define the cascade approximation operator $M=M_{m_{0}}$ as before
\begin{equation*}
Mf(x)=U^{-1}m_{0}(T)f(x)=\sqrt{2}\dsum\limits_{k\in \mathbb{Z}}a_{k}f(3x-k)%
\text{,\qquad }f\in L^{2}(\mathbb{R},(dx)^{s})\text{.}
\end{equation*}%
The condition
\begin{equation}
\frac{1}{3}\dsum\limits_{w^{3}=z}\left\vert m_{0}(w)\right\vert ^{2}=1\text{%
, a.e. }z\in \mathbb{T}\text{,}  \label{XformEq10}
\end{equation}%
will be a standing assumption on $m_{0}$. We then define the sequence $%
m_{0}^{(n)}(z):=m_{0}(z)m_{0}(z^{3})\cdots m_{0}(z^{3^{n-1}})$, and we say
that $m_{0}$ has \textit{frequency localization} if the limit of an
associated sequence of measures,
\begin{equation*}
\lim_{n\rightarrow \infty }\left\vert m_{0}^{(n)}(z)\right\vert ^{2}d\mu (z)
\end{equation*}%
exists; i.e., if there is a $\nu \in M_{1}(\mathbb{T)}$ such that
\begin{equation}
\lim_{n\rightarrow \infty }\dint\limits_{\mathbb{T}}f(z)\left\vert
m_{0}^{(n)}(z)\right\vert ^{2}d\mu (z)=\dint\limits_{\mathbb{T}}f(z)d\nu (z)
\label{XformEq11}
\end{equation}%
holds for all $f\in C(\mathbb{T)}$. Recall that if $m_{0}\in \limfunc{Lip}(%
\mathbb{T)}$ is assumed, and $R_{m_0}$ has Perron-Frobenius spectrum, then it has frequency localization, and the limit
measure $\nu $ satisfies
\begin{equation}
\lim_{n\rightarrow \infty }\mbox{dist}_{\mbox{Haus}} (d\nu
,\left\vert m_{0}^{(n)}\right\vert ^{2}d\mu )=0\text{.}
\label{XformEq12}
\end{equation}

\begin{theorem}
\label{XformTheo1}\emph{:}\textbf{\ The Dichotomy Theorem}. Let $m_{0}\in
L^{\infty }(\mathbb{T)}$ be given, and suppose \textup{(}\ref{XformEq10}%
\textup{)} holds, and
let $\nu $ be the corresponding limit measure. Assume further that, for $k\in\mathbb{Z}_+$,
\begin{equation}\label{eqextra}
\lim_{n\rightarrow\infty}\int_{\mathbb{T}}|m_0^{(n)}(z)|^2A^{(k)}_{0,0}(z)\,d\mu(z)=\nu(A^{(k)}_{0,0}),
\end{equation}
where $A$ is the matrix function defined in (\ref{XformEq6}), and

$A^{(k)}(z):=A(z)A(z^3)...A(z^{3^k})$. Let $M=M_{m_{0}}$ be the
cascade approximation operator in $L^{2}(\mathbb{R},(dx)^{s}),$
$s=\log _{3}(2)$. Then the limit
\begin{equation}
\lim_{n\rightarrow \infty }M^{n}\chi _{\mathbf{C}}\text{ exists in }L^{2}(%
\mathbb{R},(dx)^{s})  \label{XformEq13}
\end{equation}%
if and only if there is a Borel subset $E\subset \mathbb{T}$ such that $\nu
(E)=1$ (i.e., $E$ is a supporting set for $\nu $), and $m_{0}(z)=\frac{1+z^{2}%
}{\sqrt{2}}$, for all $z\in E$. In the special case where $A$ is
further assumed continuous and $m_0$ has frequency localization, the condition (\ref{eqextra}) is automatically satisfied,
\begin{equation}
Mf(x)=(M_{\mathbf{C}}f)(x)=f(3x)+f(3x-2)  \label{XformEq14}
\end{equation}%
and
\begin{equation}
M_{\mathbf{C}}\chi _{\mathbf{C}}=\chi _{\mathbf{C}}\text{.}
\label{XformEq15}
\end{equation}
\end{theorem}

\begin{proof}
Suppose first that (\ref{XformEq13}) holds; i.e., that the cascading limit
exists in $L^{2}(\mathbb{R},(dx)^{s})$. Then in particular
\begin{equation}
\lim_{n\rightarrow \infty }\left\Vert M^{n}\chi _{\mathbf{C}}-M^{n+1}\chi _{%
\mathbf{C}}\right\Vert _{L^{2}((dx)^{s})}=0\text{.}  \label{XformEq16}
\end{equation}%
We saw, using \cite{BrJo01a} that there is a measurable matrix function $A$:~%
$\mathbb{T}\rightarrow U_{3}(\mathbb{C})$ such that $m_{0}$ is the first
component in the product, matrix times vector,%
\begin{equation}
\left(
\begin{array}{c}
m_{0}(z) \\
m_{1}(z) \\
m_{2}(z)%
\end{array}%
\right) =A(z^{3})\left(
\begin{array}{c}
\frac{1+z^{2}}{\sqrt{2}} \\
z \\
\frac{1-z^{2}}{\sqrt{2}}%
\end{array}%
\right) \text{,\qquad }z\in \mathbb{T}\text{.}  \label{XformEq17}
\end{equation}%
If the functions $A_{j,k}$ denote the entries in the matrix $A$ on the
right-hand side in (\ref{XformEq17}), we saw that
\begin{equation}
p(\chi _{\mathbf{C}},M\chi _{\mathbf{C}})(z)=A_{0,0}(z)\text{,\qquad }z\in
\mathbb{T}  \label{XformEq18}
\end{equation}%
and
\begin{equation}
\left\langle M^{n}\chi _{\mathbf{C}},M^{n+1}\chi _{\mathbf{C}}\right\rangle
_{L^{2}((dx)^{s})}=\dint\limits_{\mathbb{T}}R^{n}\left( p(\chi _{\mathbf{C}%
},M\chi _{\mathbf{C}})\right) d\mu =\dint\limits_{\mathbb{T}%
}R^{n}(A_{0,0})d\mu  \label{XformEq19}
\end{equation}%
where
\begin{equation*}
R^{n}\left( A_{0,0}\right) (z)=\frac{1}{3^{n}}\dsum\limits_{w^{3^{n}}=z}%
\left\vert m_{0}^{\left( n\right) }(w)\right\vert ^{2}A_{0,0}(w)\text{%
,\qquad }z\in \mathbb{T}\text{.}
\end{equation*}%
After a change of variable in (\ref{eqextra}), we conclude from (\ref{XformEq19})
that
\begin{equation}
\lim_{n\rightarrow \infty }\left\langle M^{n}\chi _{\mathbf{C}},M^{n+1}\chi
_{\mathbf{C}}\right\rangle _{L^{2}((dx)^{s})}=\dint\limits_{\mathbb{T}%
}A_{0,0}(z)d\nu (z)=:\nu (A_{0,0})  \label{XformEq20}
\end{equation}%
and therefore
\begin{equation}
0=\lim_{n\rightarrow \infty }\left\Vert M^{n}\chi _{\mathbf{C}}-M^{n+1}\chi
_{\mathbf{C}}\right\Vert _{L^{2}((dx)^{s})}^{2}=2-2\limfunc{Re}\nu (A_{0,0})%
\text{.}  \label{XformEq21}
\end{equation}%
From (\ref{XformEq17}), we know that $\left\vert A_{0,0}\right\vert \leq 1$,
pointwise for $z\in \mathbb{T}$. From this, we get that $\nu (A_{0,0})=1$.
Hence, there is a Borel subset $E\subset \mathbb{T}$, such that $\nu (E)=1$,
and $A_{0,0}(z)=1$ for all $z\in E$. Since
\begin{equation}
\left\vert A_{0,0}(z)\right\vert ^{2}+\left\vert A_{0,1}(z)\right\vert
^{2}+\left\vert A_{0,2}(z)\right\vert ^{2}=1\text{,\qquad }z\in \mathbb{T}%
\text{,}  \label{XformEq22}
\end{equation}%
we conclude that $A_{0,1}(z)=A_{0,2}(z)=0$ for $z\in E$. Using (\ref%
{XformEq17}) again, we finally get $m_{0}(z)=\frac{1+z^{2}}{\sqrt{2}}$ for $%
z\in E$; i.e., the conclusion of the theorem holds. If $A_{0,0}$ is assumed
continuous, then $A_{0,0}=1$ on $\mathbb{T}$ since the support of $\mathbb{%
\nu }$ is all of $\mathbb{T}$. The conclusions (\ref{XformEq14})--(\ref%
{XformEq15}) in the theorem then follow.

We now turn to the converse implication: If some supporting set $E\subset
\mathbb{T}$ exists such that $A_{0,0}(z)=1$ for $z\in E$, then $\nu
(A_{0,0})=\int_{\mathbb{T}}A_{0,0}d\nu =\int_{E}A_{0,0}d\nu =\int_{E}d\nu
=\nu (E)=1$. To prove the convergence in $L^{2}((dx)^{s})$ of the cascades
in (\ref{XformEq13}), we must consider
\begin{equation}
\left\Vert M^{n}\chi _{\mathbf{C}}-M^{n+k}\chi _{\mathbf{C}}\right\Vert
_{L^{2}((dx)^{s})}^{2}=2-2\limfunc{Re}\dint\limits_{\mathbb{T}}R^{n}(p(\chi
_{\mathbf{C}},M^{k}\chi _{\mathbf{C}}))d\mu \text{,}  \label{XformEq23}
\end{equation}%
we note as in the case $k=1$, that
\begin{equation*}
p\left( \chi _{\mathbf{C}},M^{k}\chi _{\mathbf{C}}\right)
(z)=A_{0,0}^{(k)}(z)\text{.}
\end{equation*}%
If $z\in E$, then $A_{0,0}(z)=1$, and $A_{0,j}(z)=A_{j,0}(z)=0$
for $j=1,2$.

As a result, using (\ref{XformEq22}), we get
\begin{equation*}
A_{0,0}^{(2)}(z)=\dsum%
\limits_{j=0}^{2}A_{0,j}(z)A_{j,0}(z^{3})=A_{0,0}(z)A_{0,0}(z^{3})=A_{0,0}(z^{3})%
\text{,}
\end{equation*}%
and therefore
\begin{eqnarray*}
\dint\limits_{E}A_{0,0}^{(2)}(z)d\nu (z)
&=&\dint\limits_{E}A_{0,0}(z^{3})d\nu (z) \\
&=&\dint\limits_{\mathbb{T}}A_{0,0}(z^{3})d\nu (z) \\
&=&\dint\limits_{\mathbb{T}}A_{0,0}(z)d\nu (z) \\
&=&\nu (A_{0,0}) \\
&=&1\text{.}
\end{eqnarray*}

Continuing by induction, we find a supporting set, which is also denoted $E$%
, such that
\begin{equation*}
A_{0,0}^{(k)}(z)=1
\end{equation*}%
for all $z\in E$, $k=1,2,\cdots $. Since
\begin{equation*}
p(\chi _{\mathbf{C}},M^{k}\chi _{\mathbf{C}})(z)=A_{0,0}^{(k)}(z)\text{%
,\qquad }z\in \mathbb{T}\text{,}
\end{equation*}%
substitution into (\ref{XformEq23}) yields
\begin{eqnarray*}
&&\left\Vert M^{n}\chi _{\mathbf{C}}-M^{n+k}\chi _{\mathbf{C}}\right\Vert
_{L^{2}((dx)^{s})}^{2} \\
&=&2-2\limfunc{Re}\dint\limits_{\mathbb{T}}R^{n}\left( A_{0,0}^{(k)}\right)
d\mu \underset{n\rightarrow \infty }{\rightarrow }2-2\limfunc{Re}\nu \left(
A_{0,0}^{(k)}\right) \\
&=&0\text{.}
\end{eqnarray*}%
This proves convergence of the cascades, and concludes the proof of the
theorem.
\end{proof}

\section{\label{LowPass}Low pass filters}

For functions on the real line $\mathbb{R}$, and for every $N\in \mathbb{Z}%
_{+}$, $N\geq 2$, the scaling identity takes the form

\begin{equation}
\varphi (x)=D\dsum\limits_{k\in \mathbb{Z}}a_{k}\varphi (Nx-k)
\label{LowPassEq1}
\end{equation}%
where $D$ is a dimensional fixed constant. We take $D=\sqrt{N}$.
The values $a_{k}$ are called \textit{masking coefficients}, and
\begin{equation}
m_{0}(z):=D^{-1}\dsum\limits_{k\in \mathbb{Z}}a_{k}z^{k}  \label{LowPassEq2}
\end{equation}%
is the corresponding \textit{low-pass filter}. The terminology is from
graphics algorithms and signal processing, and the books \cite{Dau92} and
\cite{BrJo02} explain this connection in more detail. The function $m_{0}$
is viewed as a function on $\mathbb{T=R}/2\pi \mathbb{Z}$, or alternately as
a $2\pi $-periodic function on $\mathbb{R}$, via $z:=e^{-i\theta }$, $\theta
\in \mathbb{R}$. It turns out that the regularity properties of $m_{0}$ are
significant for the spectral theoretic properties which hold for the
operators associated with $m_{0}$, specifically \textit{the cascade
subdivision operator}, and the \textit{Ruelle transfer operator.} The
function spaces which serve as repository for the function $m_{0}$ are
measurable functions on $\mathbb{T}$, for example $L^{p}(\mathbb{T)}$, $%
1\leq p\leq \infty $, the continuous functions; i.e., $C(\mathbb{T)}$, or
the Lipschitz functions $\limfunc{Lip}(\mathbb{T)}$.

We will consider low-pass filters $m_{0}$ with the following properties:

\noindent (1) $m_{0}\in \limfunc{Lip}(\mathbb{T)}$;

\noindent (2) $m_{0}$ has a finite number of zeros;

\noindent (3) $R_{m_{0}}(1)=1$.

\begin{proposition}
\label{LowPassProp1}Suppose $m_{0}$ satisfies \textup{(1), (2), (3) and }%
\begin{equation}
\dim \left\{ h\in C(\mathbb{T)\mid }R_{m_{0}}(h)=h\right\} \geq 2\text{.}
\label{LowPassEq3}
\end{equation}%
Then there exists an $(m_{0},N)$-cycle \textup{(}i.e., a set
$\{z_{1},\ldots z_{p}\}$ with $z_{i}^{N}=z_{i+1}$,
$z_{p}^{N}=z_{1}$ and $\left\vert m_{0}(z_{i})\right\vert ^{2}=N$
for all $i$\textup{)}.
\end{proposition}

\begin{proof}
Let $h\in C(\mathbb{T)}$ satisfy $R_{m_{0}}(h)=h$, and $h$ non-constant.
Taking the real or imaginary part, we may assume that $h$ is real valued.
Also, replacing $h$ by $\left\Vert h\right\Vert _{\infty }-h$, we may assume
$h\geq 0$, and that $h$ has some zeros.

We prove that all the zeros of $h$ must be cyclic points.

Suppose not, and let $z_{0}\in \mathbb{T}$ be a zero of $h$ which is not on
a cycle. Then, $w_{1}^{N^{\ell _{1}}}=z_{0}$ and $w_{2}^{N^{\ell
_{2}}}=z_{0} $ with $\ell _{1}\neq \ell _{2}$ implies $w_{1}\neq w_{2}$.
Otherwise we have for some $\ell _{1}<\ell _{2}$, $z_{0}^{N^{\ell _{2}-\ell
_{1}}}=w_{1}^{N^{\ell _{2}}}=z_{0}$ so $z_{0}$ is a cyclic point.

We say that $w$ is at level $\ell $ if $w^{N^{\ell }}=z_{0}$. The previous
remark shows that $\ell $ is uniquely determined by $w$.

Since $h(z_{0})=0$, it follows that
\begin{equation*}
\frac{1}{N}\sum_{w^{N}=z_{0}}\left\vert m_{0}(w)\right\vert ^{2}h(w)=0\text{,%
}
\end{equation*}%
so $\left\vert m_{0}(w)\right\vert ^{2}h(w)=0$ for all $w$ with $w^{N}=z_{0}$%
. Not all such $w$'s can have $m_{0}(w)=0$ because $R_{m_{0}}(1)=1$. Thus
there is a $z_{1}$ with $z_{1}^{N}=z_{0}$ and $h(z_{1})=0$.

By induction, there is a $z_{n+1}$ with $z_{n+1}^{N}=z_{n}$ and $%
h(z_{n+1})=0 $.

Now, $m_{0}$ has only finitely many zeros, so from some level on,
there are no zeros of $m_{0}$; i.e., if $w^{N^{\ell }}=z_{0}$ with
$\ell \geq \ell _{0}$ then $m_{0}(w)\neq 0$. But then look at
those $w$'s with $w^{N}=z_{\ell
_{0}} $. Since $h(z_{\ell _{0}})=0$ and $m_0(w)\neq 0$, it follows that $%
h(w)=0 $. By induction, $h(w)=0$ for all $w$ with $w^{N^{n}}=z_{0}$ and $n$ large enough.

However, these $w$'s form a dense set because for any interval $(a,b)\subset
\mathbb{T}$ there is an $m$ big enough such that $\tau ^{m}(a,b)=\mathbb{%
T\ni }z_{\ell _{0}}$ (where $\tau (z)=z^{N}$) so there is a $w\in (a,b)$
with $w^{N^{m}}=z_{0}$. Since $h$ is continuous, it follows that $h=0$, a
contradiction.

Thus, all zeros of $h$ are cyclic points. In particular $z_{0},$ $%
z_{1},\ldots z_{n},\cdots $ are cyclic.

Since $z_{0}$ and $z_{1}$ are cyclic and $z_{1}^{N}=z_{0}$ (hence they are on
the same cycle), if $w\neq z_{1}$, $w^{N}=z_{0}$ then $w$ is not cyclic so $%
h(w)\neq 0$, and therefore $m_{0}(w)=0$. But this implies (from $%
R_{m_{0}}(1)=1$) that $\left\vert m_{0}(z_{1})\right\vert ^{2}=N$. We can do
the same for all terms of the cycle generated by $z_{1}$ and the conclusion
of the proposition follows.
\end{proof}

\begin{remark}
\label{LowPassRem1}It is known generally that the dimension of the
eigenspace in \textup{(}\ref{LowPassEq3}\textup{) }depends on the metric
properties of orbits in $\mathbb{T}$ under $z\rightarrow z^{N}$ where $N\geq
2$ is fixed. These orbits are called cycles. The solutions $h$ to
\begin{equation}
R_{m_{0}}(h)=h  \label{LowPassEq3a}
\end{equation}%
are called $\left\vert m_{0}\right\vert ^{2}$-harmonic functions and are
important in the general theory of branching processes. Their significance
for the present discussion is noted in \cite{Jorgen01}. There we prove that
each solution $h\in L^{1}(\mathbb{T})$, $h\geq 0$, $h\neq 0$, to equation
\textup{(}\ref{LowPassEq3a}\textup{)}\emph{;} i.e., a non-negative harmonic
function, is naturally associated with a system $(H,U,T,\varphi )$ where $H$
is a Hilbert space, $U$ and $T$ are unitary operators in $H$ satisfying $%
UTU^{-1}=T^{N}$\textup{(}i.e., equation \textup{(}\ref{ZakEq0}\textup{))},
and the vector $\varphi \in H$, $\varphi \neq 0$, satisfies the general
scaling identity
\begin{equation*}
U\varphi =m_{0}(T)\varphi \text{,\quad (i.e., the abstract form of (\ref%
{IntroEq3a}) or (\ref{ZakEq13}).)}
\end{equation*}%
See also \cite{BDP}.\par $T$ generates a representation of
$L^\infty(\mathbb{T})$ by
$$\pi(f)=f(T),\quad(f\in L^{\infty}(\mathbb{T})).$$
This representation satisfies the commutation relation
$$U\pi(f)U^{-1}=\pi(f(z^N)),\quad(f\in\ L^{\infty}(\mathbb{T})).$$
Iterating the scaling identity one has
$$U^n\varphi=\pi(m_0^{(n)})\varphi,\quad(n\geq 0).$$
Moreover, the system $(H,U,\pi,\varphi )$ is determined from $h$ in \textup{(}%
\ref{LowPassEq3a}\textup{)} up unitary equivalence and it is
called the wavelet representation associated to $(m_0,h)$.

Returning to the form $p(\cdot ,\cdot )$ in \textup{(}\ref{ZakEq1}\textup{)}%
, we note that $h$ is related to the new data by the two formulas,
\begin{equation*}
h=p(\varphi ,\varphi )
\end{equation*}%
and
\begin{equation*}
\dint\limits_{\mathbb{T}}\xi (z^{N})p(U\varphi ,U\varphi )(z)d\mu
(z)=\dint\limits_{\mathbb{T}}\xi (z)h(z)d\mu (z)\text{,\qquad for all }\xi
\in C(\mathbb{T)}\text{.}
\end{equation*}

The paper \cite{CoRa90} treats a general form of the problem \textup{(}\ref%
{LowPassEq3a}\textup{)}, and these authors state that the number of $%
(m_{0},N)$-cycles on $\mathbb{T}$ equals the dimension of the eigenspace in
\textup{(}\ref{LowPassEq3}\textup{)}\emph{;} i.e., the space of continuous $%
R_{m_{0}}$-harmonic functions. Recall the points $\{z_{i}\}$ on an $%
(m_{0},N) $-cycle satisfy $\left\vert m_{0}(z_{i})\right\vert =\sqrt{N}$, $%
z_{i+1}=z_{i}^{N}$, and $z_{k}=z_{0}$, if $k$ is the length of the cycle.
The following example shows that this result from \cite{CoRa90} is in need
of a slight correction$:$ Take for example $m_{0}(z)=\frac{1+z^{2}}{\sqrt{2}}
$, $N=3$, where there are no $(m_{0},3)$-cycles. What is true is that, if
the dimension in \textup{(}\ref{LowPassEq3}\textup{)} is $>$ $1$,or if there
are eigenvalues $\lambda \in \mathbb{T}\diagdown \{1\}$, then there must be $%
(m_{0},N)$-cycles, and the arguments in \cite{CoRa90} work$:$ All the
invariant measures are supported on cycles. But if the dimension in \textup{(%
}\ref{LowPassEq3}\textup{)} is $1$, then there may, or may not, be $(m_0,N)$%
-cycles. If not, then the invariant measures have support equal to $\mathbb{T%
}$. If $\{1\}$ is an $(m_0,N)$-cycle, then Dirac's $\delta _{1}$
is an invariant measure.
\end{remark}

We stress this distinction because it is central to explaining our
dichotomy; i.e., explaining when a non-zero solution $\varphi $
exists to the scaling equation (\ref{ZakEq13}). If the given
filter $m_{0}$ has an
$(m_0,N)$-cycle of length $1$, then there is a scaling function $\varphi $ in $L^{2}(\mathbb{R)}$%
, and we are in the classical (non-fractal) case. If $m_{0}$ and $N$ are
fixed, but there is no $(m_{0},N)$-cycle on $\mathbb{T}$, then $\varphi $ is
instead in one of the $\mathcal{H}^{s}$-Hilbert spaces, $0<s<1$, with $s$
depending on the chosen particular iterated function system (IFS).

\begin{proposition}
\label{LowPassProp2}Let $m_{0}$ be a filter that satisfies \textup{(1), (2),
(3)}. If there exists a $\lambda \neq 1$ with $\left\vert \lambda
\right\vert =1$ and $h\in C(\mathbb{T)}$, $h\neq 0$, such that $%
R_{m_{0}}(h)=\lambda h$, then there exists an $(m_{0},N)$-cycle.
\end{proposition}

\begin{proof}
We know that the peripheral eigenvalue spectrum of $R_{m_{0}}$ is a finite
union of cyclic subgroups of $\mathbb{T}$ (see section 4.5 in \cite{BrJo02}%
). Hence $\lambda ^{n}=1$ for some $n$.

Then $R_{m_{0}}^{n}(h)=h$ so $R_{m_{0}^{(n)}}(h)=h$ and $h$ is not a
constant. But $m_{0}^{(n)}$ is Lipschitz, it has finitely many zeros and $%
R_{m_{0}^{(n)}}(\hat{1})=\hat{1}$ (the scale for $m_{0}^{(n)}$ is $N^{n}$).

Using proposition \ref{LowPassProp1}, it follows that there is an
$(m_{0}^{(n)},N^n)$-cycle; i.e.,
there exist points $z_{1},\ldots ,z_{p}$ on $\mathbb{T}$ with z$%
_{i}^{N^{n}}=z_{i+1},z_{p}^{N^{n}}=z_{1},$ and $\left\vert
m_{0}^{(n)}(z_{i})\right\vert ^{2}=N^{n}$.

But $\left\vert m_{0}\right\vert ^{2}\leq N$ (since $R_{m_{0}}(1)=1$) and
this implies that
\begin{equation*}
\left\vert m_{0}(z_{i})\right\vert ^{2}=\left\vert
m_{0}(z_{i}^{N})\right\vert ^{2}=\ldots =\left\vert
m_{0}(z_{i}^{N^{n-1}})\right\vert ^{2}=N
\end{equation*}%
and therefore $z_{i}$ will generate an $(m_{0},N)$-cycle.
\end{proof}

\begin{theorem}
\label{LowPassTheo1}Suppose $m_{0}$ satisfies \textup{(1), (2), (3)}. Then %
exactly one of the following affirmations is true:

\noindent \textup{(i)} There exists an $(m_{0},N)$-cycle. In this
case, the invariant measures $\nu $ with
\begin{equation*}
\nu (R_{m_{0}}(f))=\nu (f),(f\in C(\mathbb{T))}
\end{equation*}%
are atomic supported on the $(m_{0},N)$-cycles. The spectrum of
$R_{m_{0}}$ is computed in \cite{BrJo02} and \cite{Dut1}. The
wavelet representation
associated to $(m_{0},1)$ is a direct sum of cyclic amplifications of $L^{2}(%
\mathbb{R)}$ (see \cite{BDP}\textup{).}

\noindent \textup{(ii)} There are no $(m_{0},N)$-cycles. In this
case there are no eigenvalues for $R_{m_{0}}\mid _{C(\mathbb{T)}}$
with $\left\vert \lambda \right\vert =1$ other then $\lambda =1;$
$1$ is a simple eigenvalue. There exists a unique probability
measure $\nu $ on $\mathbb{T}$ which is
invariant for $R_{m_{0}}$ (i.e., $\nu (R_{m_{0}}(f))=\nu (f)$ for $f\in C(%
\mathbb{T))}$.%
\begin{equation}
\lim_{n\rightarrow \infty }R_{m_{0}}^{n}(f)=\nu (f)\text{ uniformly }f\in C(%
\mathbb{T)}\text{.}  \label{LowPassEq4}
\end{equation}
\end{theorem}

\begin{proof}
When there are no $(m_{0},N)$-cycles, proposition \ref{LowPassProp1} and \ref%
{LowPassProp2} show that there are no peripheral eigenvalues other than $1$
and $1$ is a simple eigenvalue. The statements about $\nu $ follow from
theorem 3.4.4 and proposition 4.4.4 in \cite{BrJo02} and their proofs.
\end{proof}

\begin{remark}
When $m_{0}$ satisfies \textup{(1), (2), (3) and }$1$ is a simple
eigenvalue for $R_{m_{0}}\mid _{C(\mathbb{T)}}$ then the invariant
measure $\nu $ is
unique and (\ref{LowPassEq4}) holds \textup{(}see \cite{BrJo02} and \cite{Dut1}, \cite%
{BDP}\textup{). In the case of wavelet filters in }$L^{2}(\mathbb{R)}$, $%
m_{0}$ satisfies the extra condition%
\begin{equation*}
m_{0}(1)=\sqrt{N}
\end{equation*}%
so $\{1\}$ is an $(m_{0},N)$-cycle. The measure $\nu $ is simply
the Dirac measure $\delta _{1}$.
\end{remark}

\begin{lemma}
\label{LowPassLem1}Let $m_{0}$ be a filter that satisfies
\textup{(1), (2), (3)}. Assume in addition that $1$ is a simple
eigenvalue for $R_{m_{0}}\mid
_{C\left( \mathbb{T}\right) }$. Consider the wavelet representation \textup{(%
}$H,U,\pi ,\varphi $\textup{)} associated to
\textup{(}$m_{0},1$\textup{)} as in remark \ref{LowPassRem1}.

Then for all $\xi \in H$ with $\left\Vert \xi \right\Vert =1$ and all $%
f\in C\left( \mathbb{T}\right) $,
\begin{equation*}
\lim_{n\rightarrow \infty }\left\langle \xi \mid U^{-n}\pi
(f)U^{n}\xi \right\rangle =\nu (f)\text{.}
\end{equation*}
\end{lemma}

\begin{proof}
First take $\xi $ of the form $\xi =U^{-m}\pi (g)\varphi $ with $m\in
\mathbb{Z}$, $g\in C\left( \mathbb{T}\right) $. Then, for $n>m$:
\begin{eqnarray*}
&&\left\langle U^{-m}\pi
(g)\varphi \mid U^{-n}\pi (f)U^{n}U^{-m}\pi (g)\varphi  \right\rangle \\
&=&\left\langle \pi \left( g\left( z^{N^{n-m}}\right)
m_{0}^{(n-m)}(z)\right) \varphi \mid \pi \left( f(z)g\left(
z^{N^{n-m}}\right)
m_{0}^{(n-m)}(z)\right) \varphi  \right\rangle \\
&=&\dint\limits_{\mathbb{T}}f(z)\left\vert g\left( z^{N^{n-m}}\right)
\right\vert ^{2}\left\vert m_{0}^{(n-m)}(z)\right\vert ^{2}d\mu \\
&=&\dint\limits_{\mathbb{T}}\left\vert g(z)\right\vert
^{2}R_{m_{0}}^{(n-m)}f(z)d\mu \text{.}
\end{eqnarray*}%
Since $\left\Vert \xi \right\Vert =1$, it follows that $\left\Vert \pi
(g)\varphi \right\Vert =1$ so $\int_{\mathbb{T}}\left\vert g(z)\right\vert
^{2}d\mu =1$.

Also, from (\ref{LowPassEq4}), $\lim_{n\rightarrow \infty
}R_{m_{0}}^{n-m}(f)(z)=\nu (f)$ uniformly.

Therefore
\begin{equation*}
\lim_{n\rightarrow \infty }\left\langle \xi \mid U^{-n}\pi
(f)U^{n}\xi \right\rangle =\dint\limits_{\mathbb{T}}\left\vert
g(z)\right\vert ^{2}\nu (f)d\mu =\nu (f)\text{.}
\end{equation*}%
Now take $\xi \in H$ arbitrarily, $\left\Vert \xi \right\Vert =1$. We can
approximate $\xi $ by a sequence $(\xi _{j})_{j}$ of the form mentioned
before, with $\left\Vert \xi _{j}\right\Vert =1$.

Fix $\epsilon >0$. Then there is a $j$ such that $\left\Vert \xi
_{j}-\xi \right\Vert <(\frac{\epsilon }{3})\left\Vert f\right\Vert
_{\infty }$ and there is an $n_{\epsilon }$ such that, for $n\geq
n_\epsilon$,
\begin{equation*}
\left\vert \left\langle \xi_j \mid U^{-n}\pi (f)U^{n}\xi _{j}\right\rangle
-\nu (f)\right\vert <\frac{\epsilon }{3}\text{.}
\end{equation*}%
Then
\begin{gather*}
\left\vert \left\langle \xi \mid U^{-n}\pi (f)U^{n}\xi
\right\rangle -\nu
(f)\right\vert \\
\leq \left\vert \left\langle \xi \mid U^{-n}\pi (f)U^{n}\xi
\right\rangle
-\left\langle \xi  \mid U^{-n}\pi (f)U^{n}\xi _{j}\right\rangle \right\vert \\
+\left\vert \left\langle \xi \mid U^{-n}\pi (f)U^{n}\xi _{j}
\right\rangle
-\left\langle \xi_j \mid U^{-n}\pi (f)U^{n}\xi _{j} \right\rangle \right\vert \\
+\left\vert \left\langle \xi_j \mid U^{-n}\pi (f)U^{n}\xi
_{j}\mid \right\rangle -\nu (f)\right\vert \\
\leq \left\Vert f\right\Vert _{\infty }\left\Vert \xi -\xi _{j}\right\Vert
\left\Vert \xi \right\Vert +\left\Vert f\right\Vert _{\infty }\left\Vert
\xi _{j}\right\Vert \left\Vert \xi -\xi _{j}\right\Vert +\frac{\epsilon }{3}%
<\epsilon
\end{gather*}
\end{proof}

\begin{theorem}
\label{LowPassTheo2}Let $m_{0},m_{0}^{\prime }$ be two filters satisfying
\textup{(1), (2), (3)} and suppose that $1$ is a simple eigenvalue for $%
R_{m_{0}}$ and $R_{m_{0}^{\prime }}$ on $C(\mathbb{T)}$. Let $\nu $ and $\nu
^{\prime }$ be the invariant probability measures for $R_{m_{0}}$ and $%
R_{m_{0}^{\prime }}$ respectively and let \textup{(}$H,U,\pi ,\varphi $%
\textup{)}, \textup{(}$H^{\prime },U^{\prime },\pi ^{\prime },\varphi
^{\prime }$\textup{)} be the wavelet representations associated to \textup{(}%
$m_{0},1$\textup{)} and \textup{(}$m_{0}^{\prime },1$\textup{)} respectively.

If $\nu \neq \nu ^{\prime }$ then the two wavelet representations are
disjoint.
\end{theorem}

\begin{proof}
Suppose the representations are not disjoint, then there is a partial
isometry $W$ from $H$ to $H^{\prime }$, where $H$ and $H^{\prime }$ are the
respective Hilbert spaces, and $W\neq 0$. Take $\xi $ in the initial space
of $W$, $\left\Vert \xi \right\Vert =1$; then $\left\Vert W\xi \right\Vert
=1 $. Using lemma \ref{LowPassLem1} we have for all $f\in C\left( \mathbb{T}%
\right) $
\begin{eqnarray*}
\nu (f) &=&\lim_{n\rightarrow \infty }\left\langle \xi \mid
U^{-n}\pi (f)U^{n}\xi
 \right\rangle \\
&=&\lim_{n\rightarrow \infty }\left\langle W\xi \mid WU^{-n}\pi
(f)U^{n}\xi \right\rangle \\
&=&\lim_{n\rightarrow \infty }\left\langle W\xi \mid U^{\prime
-n}\pi ^{\prime
}(f)U^{\prime n}W\xi \right\rangle \\
&=&\nu ^{\prime }(f)\text{.}
\end{eqnarray*}%
Thus $\nu =\nu ^{\prime }$.
\end{proof}

\begin{corollary}
\label{LowPassCor1}Let $m_{0}$ be a filter that satisfies \textup{(1), (2), (3)}, and
suppose $1$ is a simple eigenvalue for $R_{m_{0}}$ on $C\left( \mathbb{T}%
\right) $. Let $\nu $ be the invariant measure for $R_{m_{0}}$ and let
\textup{(}$H,U,\pi ,\varphi $\textup{)} be the wavelet representation
associated to \textup{(}$m_{0},1$\textup{).} Suppose $\varphi ^{\prime }\in
H $ is another orthogonal scaling function with filter $m_{0}^{\prime }$.
Then $1$ is a simple eigenvalue for $R_{m_{0}^{\prime }}$ on $C(\mathbb{T)}$%
. If $m_{0}^{\prime }$ satisfies also \textup{(1), (2)} then the
invariant measure $\nu ^{\prime }$ for $R_{m_{0}^{^{\prime }}}$is
equal to the one for $R_{m_{0}}$, i.e., $\nu'=\nu$.
\end{corollary}

\begin{proof}
Repeating the calculation given in the proof of lemma \ref{LowPassLem1} we
have%
\begin{eqnarray*}
\nu (f) &=&\lim_{n\rightarrow \infty }\left\langle U^{-m}\pi
(g)\varphi ^{\prime } \mid U^{-n}\pi
(f)U^{n}(U^{-m}\pi (g)\varphi ^{\prime }) \right\rangle \\
&=&\lim_{n\rightarrow \infty }\dint\limits_{\mathbb{T}}\left\vert
g(z)\right\vert ^{2}R_{m_{0}^{\prime }}^{n-m}(f)d\mu
\end{eqnarray*}%
For all $f$, $g\in C\left( \mathbb{T}\right) $, $m\in \mathbb{Z}$, $\int_{%
\mathbb{T}}\left\vert g\right\vert ^{2}d\mu =1$. Suppose $h\in C\left(
\mathbb{T}\right) $, $h$ non constant with $R_{m_{0}^{\prime }}(h)=h$. Then
\begin{equation*}
\nu (h)=\dint\limits_{\mathbb{T}}\left\vert g(z)\right\vert ^{2}h(z)d\mu
\end{equation*}%
for all $g\in C\left( \mathbb{T}\right) $ with $\int_{\mathbb{T}}\left\vert
g(z)\right\vert ^{2}d\mu =1$ then $h$ is constant.

The last assertion follows directly from theorem \ref{LowPassTheo2}.
\end{proof}

\begin{corollary}
\label{LowPassCor2}Let $m_{0}$ be a filter that satisfies
\textup{(1), (2), (3) }and suppose there are no
$(m_{0},N)$-cycles. Then the wavelet representation associated to
\textup{(}$m_{0},1$\textup{)} is disjoint from the classical
wavelet representation on $L^{2}\left( \mathbb{R}\right) $.
\end{corollary}

\begin{proof}
Since $\delta _{1}$ is not invariant for $R_{m_{0}}$, everything follows
from theorem \ref{LowPassTheo2}.
\end{proof}

\textbf{Acknowledgment} Toward the end of our work on this paper
we had some very helpful discussions with Professors Steen
Pedersen, Yang Wang, Roger Nussbaum and Richard Gundy on spectral
theory of the Ruelle operators associated with IFSs and on
infinite-product formulas. We thank them for many helpful
suggestions. Corrections and helpful suggestions from the two
referees are much appreciated.

\end{document}